\documentclass[11pt]{amsart}

\usepackage{amssymb}
\usepackage{enumitem}
\usepackage{graphicx}
\usepackage{amsthm}
\usepackage{amsmath}
\usepackage{mathrsfs}
\usepackage{amsfonts}
\usepackage{mathtools}
\usepackage{bm}
\usepackage{tikz} 
\usepackage{float}
\usepackage{multirow}
\usepackage{hhline}
\usepackage{makecell}
\usepackage{caption}
\usepackage[a4paper,margin=3.4cm]{geometry}

\usetikzlibrary{matrix,decorations.pathmorphing, positioning}
\usetikzlibrary{arrows,automata,backgrounds,petri}
\usetikzlibrary{shapes.multipart}
\usetikzlibrary{decorations.pathreplacing}
\usetikzlibrary{calc, trees}
\usepackage{amsrefs}

\newcommand{\G}{\mathbf{G}}

\def\REP{\frac{2124\sqrt{2}+48\sqrt{238}}{1177}}

\makeatletter
\@namedef{subjclassname@2020}{%
  \textup{2020} Mathematics Subject Classification}
\makeatother

\newtheorem{theorem}{Theorem}[section]
\newtheorem{corollary}[theorem]{Corollary}
\newtheorem{lemma}[theorem]{Lemma}

\theoremstyle{definition}

\newtheorem{remark}[theorem]{Remark}
\newtheorem{example}[theorem]{Example}

\numberwithin{equation}{section}

\begin{document}

\title[{The Markoff and Lagrange spectra
on the Hecke group $H(4)$}]{The Markoff and Lagrange spectra \\ on the Hecke group $\mathbf H_4$} 

\author{Dong Han Kim}
\address{Department of Mathematics Education, Dongguk University - Seoul, 30 Pildong-ro 1-gil, Jung-gu, Seoul, 04620 Korea} 
\email{kim2010@dgu.ac.kr} 

\author{Deokwon Sim}
\address{Machine Learning TU, Samsung Advanced Institute of Technology,
130 Samsung-ro, Yeongtong-gu, Suwon-si, Gyeonggi-do, 16678 Korea} 
\email{deokwon.sim@snu.ac.kr}

\begin{abstract}
We consider the Markoff spectrum and the Lagrange spectrum on the Hecke group $\mathbf H_4$. 
They are identical to the Markoff and Lagrange spectra of the unit circle. 
The Markoff spectrum on $\mathbf H_4$ is also known as the Markoff spectrum of index 2 sublattices by Vulakh 
and the Markoff spectrum of 2-minimal forms or $C$-minimal forms by Schmidt. 
They characterized the spectrum up to the first accumulation point, independently.
We show that, after the first accumulation point, 
both spectra have positive Hausdorff dimension. 
Then we find gaps in the spectra and give a bound on Hall's ray.  
\end{abstract}



\maketitle

\section{Introduction}

For an irrational number $\xi$, the Lagrange number
$L(\xi)$ is defined as the supremum of all $L$ such that
$$
\left| \xi - \frac{p}{q} \right| < \frac{1}{L q^2}
$$
holds for infinitely many rational numbers $p/q$.
The classical Lagrange spectrum is the set of Lagrange numbers, i.e.,
\begin{equation}\label{def:L}
\mathscr{L}_0 : = \left\{ \limsup_{p/q \in \mathbb Q} \left( q^2 \left| \xi - \frac pq \right| \right)^{-1} \, \Big| \, \xi \in \mathbb R \setminus \mathbb Q\right\}
= \{ L(\xi) \, | \, \xi \in \mathbb R \setminus \mathbb Q\} .
\end{equation}
The Markoff spectrum is defined as the set of reciprocals of the infimum of the non-zero values of indefinite quadratic forms $f(x,y)= a x^{2}+ b xy+ c y^{2}$ with real coefficients, normalized by the square root of their discriminants $\delta(f) = b^{2}-4ac > 0$, i.e.,
\begin{equation}\label{def:M}
\mathscr{M}_0 : = \left\{ \left( \inf_{(x,y) \in \mathbb Z^2 \setminus \{ (0,0)\}} \frac{|f(x,y)|}{\sqrt{\delta(f)}} \right)^{-1} 
\, \Big| \, 
\delta(f) >0 
\right\}.
\end{equation}
It is well known \cite{To55} that $\mathscr{L}_0 \subset \mathscr{M}_0$.
The classical results of Markoff in \cite{Mar79} and \cite{Mar80} show that
$$
\mathscr{L}_0 \cap [0, 3) = \mathscr{M}_0 \cap [0, 3)
=  \left\{ \sqrt{ 9 - \frac{4}{x^2}}\, | \, x \in \mathcal{M}_0  \right \},
$$
where $\mathcal M_0 = \{ 1,2,5,13,29,34,89,169, \dots \}$ is the set of elements of positive integer triples  $(x_1, x_2, x_3)$ satisfying 
\begin{equation*}
x_1^2 + x_2^2 + x_3^2 = 3 x_1 x_2 x_3.
\end{equation*}
Therefore, the smallest accumulation point of the spectra is 3.
Moreira \cite{Mor18} showed that the two spectra have positive Hausdorff dimension right after the first accumulation point 3 
and they have full dimension starting at $\sqrt{12} - \delta$ for some $\delta >0$. 
There are gaps in $\mathscr{L}_0$ and $\mathscr{M}_0$ like $\left(\sqrt{12},\sqrt{13}\right)$, which was found by Perron \cite{Pe21}. 
Note that $\sqrt{13}$ is an isolated point on both spectra.
Eventually, there exists a half infinite interval contained in the Lagrange and Markoff spectra which is called Hall's ray \cite{Hal47}. 
Hall showed that $(6,\infty) \subset \mathscr{L}_0$ and Freiman \cite{Fre75} gave the smallest possible value $c = \frac{2 221 564 069+283 748 \sqrt{462}}{491 993 569} = 4.5278\dots$ of which $[c,\infty)$ is contained in $\mathscr{L}_0$.
For the detailed discussion of the Markoff and Lagrange spectra, see \cite{Bom07}, \cite{CF89}, \cite{LMMR21}.

The Lagrange and Markoff spectra are generalized to discrete subgroups of $\mathrm{PSL}_2(\mathbb R)$, called Fuchsian groups. 
Let $\mathbf G$ be a finitely generated Fuchsian group 
acting on the upper half plane $\mathbb H$ and its boundary $\hat{\mathbb R} = \mathbb R \cup \{ \infty\}$ via linear fractional transformation
$$
\begin{pmatrix} a & b \\ c & d \end{pmatrix} \cdot z = \frac{az+b}{cz+d}.
$$
We further assume that $\infty$ is a fixed point of a parabolic element of $\mathbf G$ and let $\mathbb Q(\mathbf G)$ be the set of orbits of $\infty$ under the action of $\mathbf G$.
For a real number $\xi$ not in $\mathbb Q(\mathbf G)$, we define the Lagrange number $L_{\mathbf G} (\xi)$ by the supremum of $L$ satisfying that 
$$
\left | \xi - M \cdot \infty \right| = \left| \xi - \frac{a}{c} \right| < \frac{1}{Lc^2}
\ \text{ for infinitely many }
M = \begin{pmatrix} a& b \\ c & d\end{pmatrix} \in \mathbf G.
$$
Since
\begin{equation*}
\left | M^{-1} \cdot \xi - M^{-1} \cdot \infty \right| = \frac{1}{c^2 \left|\xi - a/c \right|},
\end{equation*}
$L_{\mathbf G}(\xi)$ is the limit superior of $\left | M^{-1} \cdot \xi - M^{-1} \cdot \infty \right|$, which is the Euclidean diameter of the image of the geodesic from $\infty$ to $\xi$ under the action of $M^{-1} \in \mathbf G$. 
We define the Lagrange spectrum of $\mathbf G$ as
\begin{equation}\label{Lno}
\mathscr L(\mathbf G) = \left\{ L_{\mathbf G} ( \xi) \, | \, \xi \in \mathbb R \setminus \mathbb Q ( \mathbf G) \right\}.
\end{equation}

Let $f(x,y) = a x^2 + b xy + c y^2$ be an indefinite quadratic form with real coefficients.
For each quadratic form $f$, we associate a geodesic in $\mathbb H$ with end points $\xi,\eta \in \hat{\mathbb R}$, $\xi \ne \eta$ satisfying 
$$ \frac{|f(x,y)|}{\sqrt{\delta(f)}} = \frac{\left|(x-\xi y)(x - \eta y)\right|}{|\xi - \eta|}.$$
For a matrix $M = \begin{pmatrix} a & b \\ c & d \end{pmatrix} \in \mathbf G$,
we set $f(M) := f ( a, c )$
and check that 
\begin{equation}\label{eq:M1}
M \cdot \xi - M \cdot \eta = \frac{\xi - \eta}{(c\xi+d)(c\eta+d)} \quad \text{ for } \  \xi,\eta \in \hat{\mathbb R}.
\end{equation}
Therefore, we have 
$$
\frac{\sqrt{\delta(f)}}{|f(M)|}  
= \left| M^{-1} \cdot \xi - M^{-1} \cdot \eta \right|. 
$$
We define the Markoff spectrum on group $\G$ as  
\begin{align*}
\mathscr M (\G) :&= \left\{ \sup_{M \in \G} 
 \frac{\sqrt{\delta(f)}}{|f(M)|} \ \Big| \  \delta(f) > 0 \right\} \\
&= \left\{ \sup_{M \in \G} 
 \left| M^{-1} \cdot \xi - M^{-1} \cdot \eta \right| \ \Big| \ \xi, \eta \in \hat{\mathbb R}, \xi \ne \eta \right\}.
\end{align*}
The Markoff spectrum $\mathscr M(\G)$ is the set of the supremums of the Euclidean diameters of geodesics in $\mathbb H$ under the action of $\mathbf G$. 
Note that the Lagrange spectrum on group $\G$ is 
\begin{align*}
\mathscr L (\G) 
&= \left\{ \limsup_{M \in \G} 
 \left| M^{-1} \cdot \xi - M^{-1} \cdot \infty \right| \ \Big| \  \xi \in \mathbb R \setminus \mathbb Q(\mathbf G) \right\}.
\end{align*}
For the modular group $\mathrm{PSL}_2(\mathbb Z)$, we have 
$$
\mathscr M (\mathrm{PSL}_2(\mathbb Z)) = \mathscr{M}_0, \qquad
\mathscr L (\mathrm{PSL}_2(\mathbb Z)) = \mathscr{L}_0.
$$
Some closed geodesics in $\mathbb H/\mathrm{PSL}_2(\mathbb Z)$ with low heights are given in Figure~\ref{Fig:classical}.

\begin{figure}[h]
\begin{tikzpicture}[scale=3].  
	\node[below] at (0, 1) {$i$};

	
	\draw (-1/2,1.8) -- (-1/2,.866);
	\draw (1/2,1.8) -- (1/2,.866);
    \draw (.5,.866) arc (60:120:1);

    \draw[very thick] (-.5,1.118) arc (90:63.4:1.118);
    \draw[very thick] (.5,1.118) arc (90:116.6:1.118);
\end{tikzpicture}
\qquad
\begin{tikzpicture}[scale=3].  
	\node[below] at (0, 1) {$i$};

	
	\draw (-1/2,1.8) -- (-1/2,.866);
	\draw (1/2,1.8) -- (1/2,.866);
    \draw (.5,.866) arc (60:120:1);

    \draw[very thick] (-.5,1.323) arc (110.71:69.29:1.414);
    \draw[very thick] (-.5,1.323) arc (69.29:45:1.414);
    \draw[very thick] (.5,1.323) arc (110.71:135:1.414);
\end{tikzpicture}
\qquad
\begin{tikzpicture}[scale=3]. 
	\node[below] at (0, 1) {$i$};
	
	\draw (-1/2,1.8) -- (-1/2,.866);
	\draw (1/2,1.8) -- (1/2,.866);
    \draw (.5,.866) arc (60:120:1);

    \draw[very thick] (-.5,1.658) arc (106.78:73.22:1.732);
    \draw[very thick] (-.5,1.658) arc (73.22:30:1.732);
    \draw[very thick] (.5,1.658) arc (106.78:150:1.732);
\end{tikzpicture}
\caption{Several closed geodesics on the fundamental domain of the modular group on the upper half plane. 
They have maximal heights $\sqrt 5$, $2\sqrt 2$, $2\sqrt 3$ (from left to right).
}\label{Fig:classical}
\end{figure}

In this paper, we consider the Lagrange and Markoff spectra on the Hecke group $\mathbf{H}_4$ or the hyperbolic triangle group $(2,4,\infty)$.
The Hecke group $\mathbf{H}_q$ is the subgroup of $\mathrm{PSL}_2(\mathbb R)$ generated by $S = \begin{pmatrix} 0 & -1 \\ 1 & 0 \end{pmatrix}$ and $T = \begin{pmatrix} 1 & \lambda_q \\ 0 & 1 \end{pmatrix}$, where $\lambda_q = 2 \cos \frac{\pi}{q}$ and $q \ge 3$ is an integer.
The Hecke group $\mathbf{H}_q$ has the presentation
$$
\mathbf{H}_q \cong \left \langle S, T \, | \, S^2 = I, (ST)^q = I \right \rangle,
$$
where $I$ is the identity 2 by 2 matrix.
When $q = 3$, we have $\lambda_3 = 1$ and $\mathbf{H}_3$ is the modular group $\mathrm{PSL}_2(\mathbb Z)$. 
If $q = 4$, then $\lambda_4 = \sqrt 2$. 
Moreover, it is known \cite{Par77} that 
\begin{multline*}
\mathbf{H}_4 = \left\{ \begin{pmatrix} a & \sqrt 2 b \\ \sqrt 2 c & d\end{pmatrix} \, \big| \ ad-2bc = 1, \, a, b,c,d \in \mathbb Z \right\} \\
\cup \left\{ \begin{pmatrix} \sqrt 2 a & b \\ c & \sqrt 2 d\end{pmatrix} \, \big| \ 2ad-bc = 1, \, a,b,c,d \in \mathbb Z \right\}.
\end{multline*}
Therefore, we have $\mathbb Q(\mathbf H_4) = \sqrt 2 \mathbb Q$.
The Diophantine approximation on the Hecke group $\mathbf H_q$ has also been studied using the Rosen continued fraction \cite{Ro54} (see e.g. \cite{Leh85}, \cite{RS92}, \cite{Na95}, \cite{BKS00}, \cite{Na10}, \cite{BHS13}).
The three geodesics of lowest heights in $\mathbb H/\mathbf{H}_4$ are given in Figure~\ref{Fig:mod2}.

\begin{figure}[h]
\begin{tikzpicture}[scale=2.9]  
\node[below] at (0, 1) {$i$};
	
\draw (-.7071,1.7) -- (-.7071,.7071);
\draw (.7071,1.7) -- (.7071,.7071);
\draw (.7071,.7071) arc (45:135:1);

\draw[very thick] (.7071,.7071) arc (45:135:1);
\end{tikzpicture}
\
\begin{tikzpicture}[scale=2.9]  
\node[below] at (0, 1) {$i$};
	
\draw (-.7071,1.7) -- (-.7071,.7071);
\draw (.7071,1.7) -- (.7071,.7071);
\draw (.7071,.7071) arc (45:135:1);

\draw[very thick] (-.7071,1.2247) arc (90:54.74:1.2247);
\draw[very thick] (.7071,1.2247) arc (90:125.26:1.2247);
\end{tikzpicture}
\
\begin{tikzpicture}[scale=2.9] 
\node[below] at (0, 1) {$i$};
	
\draw (-.7071,1.7) -- (-.7071,.7071);
\draw (.7071,1.7) -- (.7071,.7071);
\draw (.7071,.7071) arc (45:135:1);

\draw[very thick] (-.7071,1) arc (136.69:75.96:1.4577);
\draw[very thick] (-.7071,1) arc (43.31:34.45:1.4577);
\draw[very thick] (.7071,1.414) arc (104.04:136.69:1.4577);
\draw[very thick] (-.7071,1.414) arc (75.96:43.31:1.4577);
\draw[very thick] (.7071,1) arc (136.69:145.55:1.4577);
\draw[very thick] (.7071,1) arc (43.31:104.04:1.4577);
\end{tikzpicture}
\caption{Three closed geodesics in the fundamental domain of group $\mathbf{H}_4$ on the upper half plane with lowest heights.}
\label{Fig:mod2}
\end{figure}

The minimum of the Lagrange spectrum, which is called Hurwitz's constant, for the Hecke group $\mathbf{H}_q$ was studied in \cite{Leh85} and \cite{HS86}. 
In particular, if $q$ is even, then the minimum of the Lagrange spectrum $\mathscr L(\mathbf H_q)$ is always equal to 2. 
Series \cite{Ser88} examined the discrete part of the Markoff spectrum on $\mathbf H_5$. 
The Markoff spectra on general Hecke groups were studied in \cite{Vul97}.

The discrete part of the Markoff spectrum on the Hecke group $\mathbf H_4$ has been studied  by Schmidt and Vulakh independently. 
It is known as the Markoff spectrum of 2-minimal forms by Schmidt \cite{Sch76} and the Markoff spectrum on the sublattice of index 2 by Vulakh (\cite{Vul71}; see also \cite{Mal77}).
It is also identical to the Markoff spectrum on the unit circle $S^1$,
related to the Diophantine approximation of the points $\mathbb R^2 \cap S^1$ by points of $\mathbb Q^2 \cap S^1$
(\cite{Kop85}, \cite{CK24}).
We will call $(x; y_1, y_2)$ a \emph{Vulakh-Schmidt triple} if $(x; y_1, y_2)$ is a positive integer triple satisfying
\begin{equation*}\label{MarkoffEquationIntro}
2x^2 + y_1^2 + y_2^2 = 4 x y_1 y_2.
\end{equation*}
We set 
$$\mathcal N = \{ 1, 5, 29, 65, 169, 349, \dots\} \ \text{ and } \
\mathcal M =\{ 1, 3, 11, 17, 41, 59, \dots\}
$$
as the sets of $x$'s and $y_i$'s ($i=1,2$) in the Vulakh-Schmidt triple  respectively. 
The spectral values less than $2 \sqrt 2$ are given in \cite{Sch76} (see also \cite{Mal77}) as 
$$
\mathscr M (\mathbf{H}_4) \cap \left[0,2\sqrt 2\right) = 
\left\{ \sqrt{ 8 - \frac{2}{x^2}} \ 
\Big| \ x \in \mathcal N \right \}
 \cup
\left \{ \sqrt{ 8 - \frac{4}{y^2}} \ \Big| \ 
y \in \mathcal M \right\}.
$$
Therefore, the first accumulation point of $\mathscr M(\mathbf H_4)$ is $2 \sqrt 2$.
The discrete part of the Lagrange spectrum
$\mathscr L (\mathbf{H}_4) \cap \left[0,2\sqrt 2\right)$ coincides with the discrete part of the Markoff spectrum $\mathscr M(\mathbf H_4)\cap \left[0,2\sqrt 2\right)$ (see also \cite{CK23}).
Using a method similar to the classical case, we show the first theorem.
\begin{theorem}\label{thm-1}
The Markoff spectrum $\mathscr{M}(\mathbf{H}_4)$ is closed and the Lagrange spectrum $\mathscr{L}(\mathbf{H}_4)$ is contained in $\mathscr{M}(\mathbf{H}_4)$, i.e., $\mathscr{L}(\mathbf{H}_4) \subset \mathscr{M}(\mathbf{H}_4)$.
\end{theorem}
We show that, after the first accumulation point, the Lagrange spectrum has positive Hausdorff dimension.

\begin{theorem}\label{thm0}
For any $\varepsilon >0$, we have
$$\dim_H \left( \mathscr M (\mathbf H_4) \cap \left[0, 2\sqrt 2 + \epsilon\right) \right) \ge \dim_H \left( \mathscr L (\mathbf H_4) \cap \left[0, 2\sqrt 2 + \epsilon\right) \right) > 0.$$
\end{theorem}

We call an open interval $(a,b)$ a maximal gap of the spectrum if it does not intersect the spectrum and is not a proper subset of a larger gap.
We find two maximal gaps in $\mathscr{M}(\mathbf{H}_4)$ and $\mathscr{L}(\mathbf{H}_4)$ after the first accumulation point (see Figure~\ref{fig:my_label}).

\begin{theorem}\label{thm1}
The intervals 
$\left(\frac{\sqrt{238}}{5}, \sqrt{10} \right)$ and $\left( \sqrt{10}, \frac{2124 \sqrt 2+48\sqrt{238}}{1177} \right)$
are maximal gaps in $\mathscr{M}(\mathbf{H}_4)$ and $\mathscr{L}(\mathbf{H}_4)$.
\end{theorem}

We note that $\sqrt{10}$ is an isolated point. 
The two gaps in Theorem~\ref{thm1} seem to be similar to the gaps $(\sqrt{12},\sqrt{13})$ and $\left(\sqrt{13},\frac{1}{22}(9\sqrt{3}+65)\right)$ in the classical Markoff and Lagrange spectra $\mathscr{M}_0$ and $\mathscr{L}_0$ \cite{CF89}*{Lemmas 7 and 9}.

After a certain point the Lagrange spectrum $\mathscr{L}(\mathbf{H}_4)$ contains a half line, so does $\mathscr M(\mathbf{H}_4)$, 
which is called Hall's ray (see Figure~\ref{fig:my_label}). 
The existence of Hall's ray in $\mathscr{L}(\mathbf{H}_4)$ is established \cite{AMU20} in general groups.
We give a bound of the Hall's ray as follows.

\begin{figure}[h]
\begin{tikzpicture}[scale=6]
\draw[->] (1.8,0) -- (4.0,0);
\filldraw [black] (2,0) circle (0.2pt) node[below] {$2\sqrt 2$};
\filldraw [black] (1.91,0) circle (0.2pt);
\filldraw [black] (1.96,0) circle (0.2pt);
\filldraw [black] (1.99,0) circle (0.2pt);
\filldraw [black] (2.2,0) circle (0.2pt) node[below] {$\frac{\sqrt{238}}{5}$};
\filldraw (2.38,0) node[below] {$\sqrt{10}$}; \filldraw [black] (2.4,0) circle (0.2pt);
\filldraw [black] (2.5,0) circle (0.2pt); 
\filldraw (2.65,0) node[below] {$\frac{2124\sqrt{2}+48\sqrt{238}}{1177}$};
\filldraw [black] (3.7,0) circle (0.2pt) node[below] {$4\sqrt 2$};
        
\draw [decorate,decoration={brace,amplitude=3pt},xshift=0pt,yshift=0.1pt] (2.2,0) -- (2.4,0) node[black,midway,yshift=6pt]{gap};
\draw [decorate,decoration={brace,amplitude=2pt},xshift=0pt,yshift=0.1pt] (2.4,0) -- (2.5,0) node[black,midway,yshift=6pt]{gap};
\draw[very thick,black] (3.5,0) -- (4.0,0) node[black,midway,yshift=5pt] {ray};
    
\end{tikzpicture}
\caption{Gaps and a ray in $\mathscr{M}(\mathbf{H}_4)$}
\label{fig:my_label}
\end{figure}

\begin{theorem}\label{thm2}
The Lagrange spectrum $\mathscr L(\mathbf{H}_4)$ contains every real number greater than or equal to $4\sqrt2$, i.e. $[4\sqrt2,\infty)\subset \mathscr L(\mathbf{H}_4) \subset \mathscr{M}(\mathbf{H}_4)$.
\end{theorem}

In Section~\ref{Sec:Perron}, we introduce a symbolic coding for a geodesic and its endpoints in $\hat{\mathbb R}$ 
using the Hecke group $\mathbf H_4$.
A geodesic in the hyperbolic space is determined by a doubly infinite expansion and we deduce the formula of the spectral value of the Markoff and Lagrange spectra by the doubly infinite expansion of the geodesic.
We then prove Theorem~\ref{thm-1}, Theorem~\ref{thm0}, Theorem~\ref{thm1} and Theorem~\ref{thm2} in Section~\ref{Sec:closedness}, Section~\ref{Sec:dim},  Section~\ref{Sec:Gaps}, and Section~\ref{Sec:Halslray} respectively.

\section{Symbolic coding of a geodesic and the Perron formula}\label{Sec:Perron}

In this section, we introduce a symbolic coding for a geodesic and its endpoints using the Hecke group $\mathbf H_4$, following the work of Hass and Series \cite{HS86} and Series \cite{Ser88}.
We then derive the Perron formula (Theorems \ref{P1} and \ref{P2}) for $\mathbf H_4$ using this expansion. 
The $\mathbf H_4$-expansion is closely related with 
the digit expansion on the unit circle introduced by Romik \cite{Rom08},
which is also related with the even integer continued fraction or continued fractions of specific parities (see \cite{KLL22}, \cite{KLL25}).
For the connection between $\mathbf H_4$-expansion and the Rosen continued fraction, consult \cite{BKL}.

Let 
$$
T = \begin{pmatrix} 1 & \sqrt 2 \\ 0 & 1 \end{pmatrix}, 
\quad 
S = \begin{pmatrix} 0 & -1 \\ 1 & 0 \end{pmatrix}, 
\quad 
K = ST^{-1} = \begin{pmatrix} 0 & 1 \\ -1 & \sqrt 2 \end{pmatrix}.
$$
Note that $K^4 = I$.
We consider a fundamental domain $\Omega$ for $\mathbf H_4$ surrounded by geodesics given by $x=0$, $x = \sqrt 2$, $|z| = 1$ and $|z - \sqrt 2| = 1$ (Figure~\ref{Fig:P} (left)).
Let $\bm\delta_0$ be the geodesic given by the imaginary axis
and $\bm\delta_d = K^d (\bm\delta_0)$ for $d = 1,2,3$.
Let $\Delta = \Omega \cup K(\Omega) \cup K^2(\Omega) \cup K^3 (\Omega)$ be the ideal quadrilateral with edges $\bm\delta_d$ for $d = 0,1,2,3$  (Figure~\ref{Fig:P} (right)).
Let $\bm\Gamma_4$ be the subgroup of $\mathbf H_4$ generated by $K^d S K^{-d}$, $d = 0,1,2,3$.
Then $\Delta$ is a fundamental domain of $\bm\Gamma_4$ (see \cite{HS86}).

\begin{figure}[h]
\begin{center}
\begin{tikzpicture}[scale=3.6]
\node[below] at (0, 0) {$0$};
\node[below] at ({1/sqrt(2)}, 0) {$\frac 1{\sqrt 2}$};
\node[left] at (0,1) {$i$};
\node[below] at ({sqrt(2)}, 0) {$\frac{\sqrt 2}1$};


\draw[thick] (-.1, 0) -- ({sqrt(2)+.1}, 0);
\draw[thick] (0,1) -- (0,1.4);
\draw[thick] ({sqrt(2)},1) -- ({sqrt(2)},1.4);


\draw[thick,dashed,red] ({1/sqrt(2)},{1/sqrt(2)}) -- ({1/sqrt(2)},1.4);

\draw[thick] (0,1) arc (90:45:1);
\draw[thick] ({sqrt(2)},1) arc (90:135.53:1);
\end{tikzpicture}
\quad
\begin{tikzpicture}[scale=3.6]
\node[below] at (0, 0) {$0$};
\node[below] at ({1/sqrt(2)}, 0) {$\frac 1{\sqrt 2}$};
\node[left] at (0,1) {$i$};
\node[below] at ({sqrt(2)}, 0) {$\frac{\sqrt 2}1$};

\node[left] at (0, .75) {$\bm\delta_0$};
\node[below] at ({1/(2*sqrt(2))},{1/(2*sqrt(2))}) {$\bm\delta_1$};
\node[below] at ({3/(2*sqrt(2))},{1/(2*sqrt(2))}) {$\bm\delta_2$};
\node[right] at ({sqrt(2)}, .75) {$\bm\delta_3$};

\draw[thick] (-.1, 0) -- ({sqrt(2)+.1}, 0);
\draw[thick] (0,0) -- (0,1.4);
\draw[thick] ({sqrt(2)},0) -- ({sqrt(2)},1.4);

\draw[thick] (0,0) arc (180:0:{1/(2*sqrt(2))});

\draw[thick] ({sqrt(2)},0) arc (0:180:{1/(2*sqrt(2))});

\draw[thick,dashed,red] ({1/sqrt(2)},0) -- ({1/sqrt(2)},1.4);
\draw[thick,dashed,red] (0,0) arc (180:0:{1/sqrt(2)});

\draw[thick,dashed,blue] (0,1) arc (90:19.47:1);
\draw[thick,dashed,blue] ({sqrt(2)},1) arc (90:160.53:1);
\end{tikzpicture}
\end{center}

\caption{A fundamental domain $\Omega$ 
(left) and the ideal quadrilateral $\Delta$ (right)}\label{Fig:P}
\end{figure}

Let $\bm\gamma$ be an oriented geodesic with end points $\bm\gamma^-, \bm\gamma^+ \in \hat{\mathbb R}$.
We assume that neither $\bm\gamma^-$ nor $\bm\gamma ^+$ belongs to $\mathbb Q(\mathbf H_4)$.
Let $\mathscr T 
= \cup_{G \in \mathbf H_4} G(\bm\delta_0)$. 
Then, by cutting $\mathscr T$, the oriented geodesic $\bm\gamma$ is divided into geodesic segments $\dots, \bm\gamma_{-2}, \bm\gamma_{-1}, \bm\gamma_{0}, \bm\gamma_{1}, \bm\gamma_{2}, \dots$ along the orientation. 
Let $\bm\gamma^-_n, \bm\gamma^+_n \in \mathscr T$, be the two end points of the geodesic segment $\bm\gamma_n$ along the orientation of $\bm\gamma$. 
For each $n \in \mathbb Z$, there exists $M_n \in \bm\Gamma_4$ such that $\bm\gamma_n$ belongs to $M_n(\Delta)$.  
Let $e_n \in \{ 0,1,2,3\}$ be such that $\bm\gamma_n^- \in M_n (\bm\delta_{e_n})$
and define $G_n = M_n K^{e_n}$. 
Then  
$G_n \in \mathbf H_4$ and 
\begin{equation}\label{geodesic}
\bm\gamma_n^- \in G_n (\bm\delta_0) \quad \text{ and } \quad
\bm\gamma_n^+ \in G_n (\bm\delta_{d_n}) \quad  \text{ for some } \ d_n \in \{ 1,2,3\}.
\end{equation}
Since 
$$\bm\gamma^-_{n+1} = \bm\gamma_{n}^+ \in G_n( \bm\delta_{d_n} )
= G_n K^{d_n} (\bm\delta_0) = G_n K^{d_n} S (\bm\delta_0),
$$ 
we deduce that for all $n \in \mathbb Z$
\begin{equation}\label{gNn}
G_{n+1} = G_n K^{d_n} S = G_{n} N_{d_n},
\end{equation}
where we denote 
$$N_d := K^dS \quad \text{ for } d = 1,2,3.$$
Note that 
\begin{equation}\label{def:As}
N_1 = \begin{pmatrix} 1 & 0 \\ \sqrt 2 & 1 \end{pmatrix},
\quad 
N_2 = \begin{pmatrix} \sqrt 2 & 1 \\ 1 & \sqrt 2 \end{pmatrix},
\quad 
N_3 =\begin{pmatrix} 1 & \sqrt 2 \\ 0 & 1 \end{pmatrix}.
\end{equation}

For each oriented geodesic $\bm\gamma$ on $\mathbb H$, we define a two-sided infinite sequence
$(d_n)_{n\in\mathbb Z} \in \{ 1,2,3\}^{\mathbb Z}$.
We give an equivalence relation $(a_n) \sim (b_n)$ in $\{1,2,3\}^{\mathbb{Z}}$ if and only if there exists some $m \in \mathbb Z$ such that $a_{n+m} = b_n$ for all $n \in \mathbb Z$. 
Then an equivalence class of $\{1,2,3\}^{\mathbb{Z}}$ under the equivalence relation is called a \emph{doubly-infinite $\mathbf H_4$-sequence}.
A \emph{section} of a doubly-infinite $\mathbf H_4$-sequence is an element $(d_n)_{n\in\mathbb Z} \in \{ 1,2,3\}^{\mathbb Z}$ in the equivalence class. 
For each oriented geodesic $\bm\gamma$ on $\mathbb H$, we associate a doubly-infinite $\mathbf H_4$-sequence.
Figure~\ref{Fig:F} shows an example of oriented geodesic $\bm\gamma$ with $G_1 = I$ and
$$
\dots, \ d_{-1} = 1, \ d_0 = 2, \ d_1 = 3, \ d_2 = 1, \ d_3 = 3, \ \dots
$$
thus, 
$$G_{-1} = N_2^{-1}N_1^{-1},  G_{0} = N_2^{-1}, G_{1} = I,  G_{2} = N_3,  G_{3} = N_3N_1, G_{4} = N_3N_1N_3.$$

\begin{figure}[h]
\centering
\begin{tikzpicture}[scale=2.25]
\node[below] at (0, 0) {$0$};
\node[below] at (1, 0) {$\frac 1{\sqrt 2}$};
\node[below] at (2, 0) {$\frac{\sqrt 2}1$};

\node[blue] at (-1, 1.3) {$G_0(\Delta)$};
\node at (1, 1.3) {$G_1(\Delta) = \Delta$};
\node[blue] at (3, 1.3) {$G_2(\Delta)$};
\node[blue] at (2.6, .3) {$G_3(\Delta)$};
\node at (-1.45, .3) {$G_{-1}(\Delta)$};
\node at (3.3, .1) {$G_4(\Delta)$};

\node[below,red] at (2.9, 0) {$\bm\gamma^+$};
\node[below,red] at (-1.1, 0) {$\bm\gamma^-$};
\node[red] at (1,1.85) {$\bm\gamma_1$};
\node[red] at (2.8,1) {$\bm\gamma_2$};
\node[red] at (-1,1) {$\bm\gamma_0$};

\draw[thick] (-2, 0) -- (4, 0);
\draw[blue] (-2,0) -- (-2,2.1);
\draw[thick] (0,0) -- (0,2.1);
\draw[thick] (2,0) -- (2,2.1);
\draw[blue] (4,0) -- (4,2.1);

\draw[thick] (-2,0) arc (180:0:1/2);
\draw[thick] (-2,0) arc (180:0:1/4);

\draw[blue] (0,0) arc (0:180:1/2);
\draw[thick] (0,0) arc (180:0:1/2);

\draw[thick] (2,0) arc (0:180:1/2);
\draw[blue] (2,0) arc (180:0:1/2);
\draw[blue] (2,0) arc (180:0:1/4);

\draw[blue] (4,0) arc (0:180:1/2);

\draw[thick] (-1,0) arc (0:180:1/6);
\draw[thick] (3,0) arc (0:180:1/6);

\draw[thick] (-3/2,0) arc (180:0:1/12);
\draw[blue] (5/2,0) arc (180:0:1/12);

\draw[thick] (3,0) arc (0:180:1/10);
\draw[thick] (11/4,0) arc (0:180:1/24);
\draw[thick] (11/4,0) arc (180:0:1/40);


\draw[thick, red] (-1.1,0) arc (180:0:2.0);

\end{tikzpicture}
\caption{An oriented geodesic $\bm\gamma$ with 
a sequence of geodesic segments $\bm\gamma_n$.
}\label{Fig:F}
\end{figure}

From \eqref{geodesic}, we deduce that for each $n \in \mathbb Z$ the oriented geodesic $G_{n}^{-1}(\bm\gamma)$ intersects the imaginary axis $\bm\delta_0$
and satisfies  
\begin{equation}\label{endpoints}
G_{n}^{-1}(\bm\gamma^-) \in (-\infty,0) \quad \text{ and } \quad G_{n}^{-1}(\bm\gamma^+) \in (0, \infty).
\end{equation}
Suppose that $G_1 = I$. 
Then 
we have $\bm\gamma^+ \in (0, \infty)$
and by \eqref{gNn}, 
we obtain $G_{n+1} = N_{d_{1}} N_{d_{2}} \cdots N_{d_{n}} $ for $n \ge 0$.
Therefore, \eqref{endpoints} implies that for all $n \ge 1$,
$$
\bm\gamma^+ \in N_{d_{1}} N_{d_{2}} \cdots N_{d_{n}} \cdot (0,\infty).
$$

Using the symbolic coding of the geodesic, 
we have an expansion of a positive real number by the one-sided infinite sequence $(d_{n})_{n \in \mathbb N}$.
Let 
$f : [0,\infty] \to [0,\infty]$ be the map given by
$$
f(x) = \begin{cases}
N_1^{-1} \cdot x, & \text{ if } x \in \big[0,\frac{1}{\sqrt 2}\big] = N_1 \cdot [0,\infty], \\
N_2^{-1} \cdot x, & \text{ if } x \in \big[ \frac{1}{\sqrt 2}, \sqrt 2\big]  = N_2 \cdot [0,\infty], \\
N_3^{-1} \cdot x, & \text{ if } x \in \big[\sqrt 2, \infty \big]  = N_3 \cdot [0,\infty].
\end{cases}
$$
For a real number $\alpha \in [0,\infty]$,
there exists an infinite sequence $(d_n)_{n \in \mathbb N}$ satisfying 
$$
f^{n-1}(x) \in N_{d_n} \cdot [0,\infty] \quad \text{ for all }\ n \ge 1,
$$
thus
$$
x \in N_{d_1} N_{d_2} \cdots N_{d_n} \cdot [0,\infty] \quad \text{ for all }\ n \ge 1.
$$
We define the \emph{$\mathbf H_4$-expansion} of $\alpha$ as 
$$ \alpha = [d_1, d_2, d_3, \dots].$$
\begin{remark}
For $q=3$ we have $T = \begin{pmatrix} 1 & 1 \\ 0 & 1 \end{pmatrix}$ and $\mathbf H_3 = \mathrm{PSL}_2(\mathbb Z)$. 
In this case, we have 
$$
N_1 = \begin{pmatrix} 1 & 0 \\ 1 & 1 \end{pmatrix}, \quad
N_2 = \begin{pmatrix} 1 & 1 \\ 0 & 1 \end{pmatrix}.
$$
Thus, 
the $\mathbf H_3$-expansion of $\alpha$ is  
$$
\alpha = [ \overbrace{2, \dots , 2}^{a_0}, \overbrace{1, \dots 1}^{a_1}, \overbrace{2, \dots, 2}^{a_2}, \dots ] \quad \text{ for } \ \alpha = a_0 + \cfrac{1}{a_1 + \cfrac{1}{a_2+\ddots}}.
$$
\end{remark}

By the \emph{infinite $\mathbf H_4$-sequence}, we mean an element of $\{1,2,3\}^{\mathbb{N}}$.
For an infinite $\mathbf H_4$-sequence $P=(d_n)_{n\ge1}$, we write $[P] = [d_1, d_2, \dots]$.
For $d_1, \dots, d_k \in \{1,2,3\}$, we define the \emph{cylinder set} 
$$[d_1, d_2, \dots, d_k ] := N_{d_1} \cdots N_{d_k} \cdot [0,\infty].$$
For $\alpha = [d_1, d_2, \dots ]$, we have 
$\alpha \in [d_1, d_2, \dots, d_n]$  for all $n \ge 1.$
Some cylinder sets of the $\mathbf H_4$-expansion are given in Figure~\ref{Fig:cylinderset}.
We note that for each $k \ge 1$
$$
[d_1, d_2, \dots] = N_{d_1} \cdots N_{d_k} \cdot [d_{k+1}, d_{k+2}, \dots].  
$$
In particular we check 
\begin{equation}\label{digit_relation}
[1,P]=N_1 \cdot [P], \qquad 
[2,P]=N_2 \cdot [P], \qquad 
[3,P]=N_3 \cdot [P]. 
\end{equation}
and deduce that 
$$0 \le [1,P] \le \frac{1}{\sqrt 2}, \qquad
\frac{1}{\sqrt 2} \le [2,P] \le \sqrt 2, \qquad
\sqrt 2 \le [3,P]
$$
for $P\in\{1,2,3\}^{\mathbb{N}}$.

\begin{figure}[h]
\begin{tikzpicture}[xscale=3.4]
\draw[->] (1.1, 0) -- (5, 0);
\draw (1.1, 0.1) -- +(0, -0.2) node[below] {$0$};
\draw (3/2, 0.1) -- +(0, -0.2) node[below] {$\frac1{2\sqrt 2}$};
\draw (5/3, 0.1) -- +(0, -0.2) node[below] {$\frac{\sqrt2}3$};
\draw (2, 0.1) -- +(0, -0.2) node[below] {$\frac{1}{\sqrt 2}$};
\draw (7/3, 0.1) -- +(0, -0.2) node[below] {$\frac{2\sqrt2}3$};
\draw (5/2, 0.1) -- +(0, -0.2) node[below] {$\frac3{2\sqrt 2}$};
\draw (3, 0.1) -- +(0, -0.2) node[below] {$\frac{\sqrt 2}{1}$};
\draw (3.8, 0.1) -- +(0, -0.2) node[below] {$\frac{3}{\sqrt 2}$};
\draw (4.6, 0.1) -- +(0, -0.2) node[below] {$\frac{2\sqrt 2}{1}$};

\draw[<->, dotted] (1.1, 0.5) -- (3/2, 0.5) node[midway, above] {$[1, 1]$};
\draw[<->, dotted] (3/2, 0.5) -- (5/3, 0.5) node[midway, above] {$[1, 2]$};
\draw[<->, dotted] (5/3, 0.5) -- (2, 0.5) node[midway, above] {$[1, 3]$};
\draw[<->, dotted] (2, 0.5) -- (7/3, 0.5) node[midway, above] {$[2, 1]$};
\draw[<->, dotted] (7/3, 0.5) -- (5/2, 0.5) node[midway, above] {$[2, 2]$};
\draw[<->, dotted] (5/2, 0.5) -- (3, 0.5) node[midway, above] {$[2, 3]$};
\draw[<-, dotted] (3, 0.5) -- (3.8, 0.5) node[midway, above] {$[3,1]$};
\draw[<-, dotted] (3.8, 0.5) -- (4.6, 0.5) node[midway, above] {$[3,2]$};
\draw[<-, dotted] (4.6, 0.5) -- (5, 0.5) node[midway, above] {$[3,3]$};
\end{tikzpicture}
\caption{Cylinder sets on $\mathbb R$}\label{Fig:cylinderset}
\end{figure}

Since $\mathbf H_4$ is generated by $S$ and $K$, any $M \in \mathbf H_4$ takes one of the following forms 
$$
N_{d_1} \cdots N_{d_k} \ \text{ or } \ N_{d_1} \cdots N_{d_k}S
\ \text{ or } \ S N_{d_1} \cdots N_{d_k} \ \text{ or } \ SN_{d_1} \cdots N_{d_k}S.
$$
Therefore, we have $\alpha \in [0,\infty]$ belongs to $\mathbb Q(\mathbf H_4)$ if and only if
$$ \alpha = N_{d_1} \cdots N_{d_k} \cdot 0 \quad \text{ or } \quad  \alpha = N_{d_1} \cdots N_{d_k} \cdot \infty, $$
 which is equivalent to $\alpha$ being a boundary point of a cylinder set $[d_1, \dots, d_k]$.
If $\alpha$ belongs to $\mathbb R \setminus \mathbb Q(\mathbf H_4)$, then it has a unique $\mathbf H_4$-expansion
$[d_1, d_2, \dots ]$.
For the boundary points of the cylinder set, we have 
\[
0 = [1, 1 ,1, \dots] =: [1^{\infty}], \qquad \infty = [3, 3 ,3, \dots] =: [3^{\infty}].
\]
and 
$$
[d_1, d_2, \dots, d_k ] = \big[ \, [d_1, d_2, \dots, d_k, 1^\infty], \, [d_1, d_2, \dots, d_k, 3^\infty] \, \big].
$$
Therefore, if $\alpha$ belongs to $\mathbb Q(\mathbf H_4)$, then there exist up to two expressions of $\alpha$. 
For instance,
$$
\frac{1}{\sqrt 2} = [ 1, 3^\infty ] = [2, 1^\infty], \qquad 
\sqrt 2 = [ 2, 3^\infty ] = [3, 1^\infty].
$$
\begin{example}\label{exam1}
Since $[2^\infty] = N_2 \cdot [2^\infty]$, 
$[(1, 3)^\infty] = N_1 N_3 \cdot [(1, 3)^\infty]$, we have 
\begin{equation}\label{period_e1}
[2^\infty]=1, \quad  
[(1, 3)^\infty] = \frac{\sqrt 3 - 1}{\sqrt 2}.
\end{equation}
Similarly, we check 
\begin{equation}\label{period_e2}
[(1,2,3)^\infty] = \frac{\sqrt{17}-2\sqrt 2}{3}, \quad
[ (1,1,2)^\infty ] = \frac{1}{\sqrt 7 + \sqrt 2}.
\end{equation}
\end{example}

For infinite $\mathbf H_4$-sequences $P=(a_n)_{n\ge1}$ and $Q=(b_n)_{n\ge1}$, we define a combined two-sided sequence 
$$P^{*}|Q := (c_n)_{n \in \mathbb Z}, \qquad c_n = \begin{cases} b_n, & \text{ if } n \ge 1, \\ a_{-n+1}, &\text{ if } n \le 0, \end{cases}$$ 
which is an element of $\{1,2,3\}^{\mathbb{Z}}$.
Let
$$d^{\vee}= 
\begin{cases}
3 &\text{ if } \ d=1,\\
2 &\text{ if } \ d=2,\\
1 &\text{ if } \ d=3.
\end{cases}
$$
Then we have identities 
\begin{equation}\label{invN}
N_{d}^{-1}  = S N_{d^\vee} S \quad \text{ and } \quad  N_{d^\vee}  = J N_{d} J \quad \text{ where } \ J = \begin{pmatrix} 0 & 1 \\ 1 & 0 \end{pmatrix}.
\end{equation}
For a given infinite $\mathbf H_4$-sequence $P =  (a_n)_{n \ge 1}$, let $P^{\vee} = \left( a_n^\vee \right)_{n\ge 1}$. 
For a doubly-infinite $\mathbf H_4$-sequence $U$ with a section $P^*|Q$, we define $U^{\vee}$ and $U^*$ as the doubly-infinite $\mathbf H_4$-sequences with a section $(P^\vee)^*| Q^\vee$ and $Q^*|P$ respectively. 
Using \eqref{invN}, we have 
\begin{align*}
[d^\vee_1, \dots, d^\vee_k ] 
&=  J N_{d_1} \cdots N_{d_k} J \cdot [0,\infty] \\
&= \left[ \, \frac{1}{[d_1, d_2, \dots, d_k, 3^\infty]}, \, \frac{1}{[d_1, d_2, \dots, d_k, 1^\infty]}  \, \right]
\end{align*}
and 
$$ [P^\vee ] =  \frac{1}{[P]}.
$$
For an example, from \eqref{period_e1}, we have 
\begin{equation*}
[(3, 1)^\infty ] = \frac{1}{[(1, 3)^\infty ]} = \frac{\sqrt 2}{\sqrt 3 -1} = \frac{\sqrt 3 + 1}{\sqrt 2}.
\end{equation*}
We also note that 
$$
[(3, 1)^\infty ] 
= N_3 \cdot [(1, 3)^\infty] = \sqrt 2 + [(1, 3)^\infty] = \frac{\sqrt 3 + 1}{\sqrt 2}.
$$


\begin{lemma}\label{lem:std}
Suppose that $\bm\gamma$ is an oriented geodesic on $\mathbb H$ with an associated doubly infinite $\mathbf H_4$-sequence $U$.

\noindent
(i) If $\bm\gamma^- \in (-\infty, 0)$ and $\bm\gamma^+ \in (0,\infty)$, then there exists a section $P^*|Q$ of $U$ with $P = (a_n)_{n\in \mathbb N}$, $Q = (b_n)_{n \in \mathbb N}$ satisfying 
$$\bm\gamma^- = -[P] = - [a_1, a_2, \dots] \quad \text{  and } \quad \bm\gamma^+ = [Q] = [b_1, b_2, \dots]. $$
\noindent
(ii) If $\bm\gamma^- \in (0, \infty)$ and $\bm\gamma^+ \in (-\infty,0)$,
then there exists a section $P^*|Q$ of $U$ such that  
$$\bm\gamma^+ = -[P^\vee] = - [a^\vee_1, a^\vee_2, \dots] \quad \text{  and } \quad \bm\gamma^- = [Q^\vee] = [b^\vee_1, b^\vee_2, \dots]. $$
\end{lemma}

\begin{proof}
We first assume that $\bm\gamma^- \in (-\infty, 0)$, $\bm\gamma^+ \in (0, \infty)$. 
Then we choose geodesic segments $\bm\gamma_0$ and $\bm\gamma_1$ in $S(\Delta)$ and $\Delta$ respectively,
thus, $\bm\gamma_1^-\in \bm\delta_0$ and $G_1 = I$. 
By \eqref{gNn} and \eqref{invN} we obtain that  
\begin{equation*}
G_n = \begin{cases} N_{d_1} N_{d_{2}} \cdots N_{d_{n-1}}, & \text{ if } \ n \ge 2, \\  
S N_{d^\vee_{0}} N_{d^\vee_{-1}} \cdots N_{d^\vee_{n}} S, & \text{ if } \ n \le 0.
\end{cases}
\end{equation*}
Therefore, by \eqref{endpoints}, we have for all $m \ge 1$
\begin{equation*}
\bm\gamma^+ \in N_{d_1}N_{d_2} \cdots N_{d_{m}} \cdot (0,\infty), 
\quad
S \cdot \bm\gamma^- \in N_{d^\vee_0} N_{d^\vee_{-1}} \cdots N_{d^\vee_{-m}} \cdot ( 0, \infty),
\end{equation*}
which yields that
$$
\bm\gamma^+ = [d_1, d_2, d_3, \dots], \quad
\bm\gamma^- = - \frac{1}{[d^\vee_0, d^\vee_{-1}, d^\vee_{-2}, \dots]} = - [d_0, d_{-1}, d_{-2}, \dots].
$$

Next, we consider the case of $\bm\gamma^- \in (0, \infty)$, $\bm\gamma^+ \in (-\infty, 0)$.
We choose geodesic segments $\bm\gamma_0$ and $\bm\gamma_1$ in $\Delta$ and $S(\Delta)$ respectively,
thus, $\bm\gamma_1^-\in S(\bm\delta_0) = \bm\delta_0$ and $G_1 = S$. 
By \eqref{gNn} and \eqref{invN} we get  
\begin{equation*}
G_n = \begin{cases} 
S N_{d_1} N_{d_{2}} \cdots N_{d_{n-1}}, 
& \text{ if } \ n \ge 2, \\  
N_{d^\vee_{0}} N_{d^\vee_{-1}} \cdots N_{d^\vee_{n}} S, & \text{ if } \ n \le 0.
\end{cases}
\end{equation*}
Therefore, by \eqref{endpoints}, we have for all $m \ge 1$ 
\begin{equation*}
S \cdot \bm\gamma^+ \in N_{d_1}N_{d_2} \cdots N_{d_{m}} \cdot (0,\infty), \quad
\bm\gamma^- 
\in N_{d^\vee_{0}} N_{d^\vee_{-1}} \cdots N_{d^\vee_{-m}} \cdot (0, \infty),
\end{equation*}
which implies that 
\[
\bm\gamma^+ =- \frac{1}{[d_1, d_2, d_3, \dots]}  = - [d^\vee_1, d^\vee_{2}, d^\vee_{3}, \dots], \quad
\bm\gamma^- = [d^\vee_0, d^\vee_{-1}, d^\vee_{-2}, \dots]. \qedhere
\]
\end{proof}

\begin{lemma}\label{lem:min2}
Let $\bm\gamma$ be an oriented geodesic on $\mathbb H$ with two end points $\bm\gamma^-$, $\bm\gamma^+$
and let $U$ be the doubly infinite $\mathbf H_4$-sequence associated to $\bm\gamma$.
If  $\left| M \cdot \bm\gamma^+ - M \cdot \bm\gamma^- \right| > \sqrt 2$ for some $M \in \mathbf H_4$, then there exists a section $P^*|Q$ of $U$ such that
$$\left| M \cdot \bm\gamma^+ - M \cdot \bm\gamma^- \right| 
= [P] + [Q] \ \text{ or } \  [P^\vee] + [Q^\vee].
$$
\end{lemma}

\begin{proof}
Suppose that 
$$
\left| M \cdot \bm\gamma^+ - M \cdot \bm\gamma^- \right| > \sqrt 2.
$$
By replacing $M' = T^m M$ for some $m \in \mathbb Z$, we may assume that 
$$
M' \cdot \bm\gamma^- < 0,  \ M' \cdot \bm\gamma^+ > 0 \quad \text{ or } \quad  M' \cdot \bm\gamma^- > 0,  \ M' \cdot \bm\gamma^+ < 0.
$$
Let $\tilde{\bm\gamma} = M'(\bm\gamma)$. 
Then, $U$ is also the associated doubly infinite $\mathbf H_4$-sequence of $\tilde{\bm\gamma}$.
By Lemma~\ref{lem:std}, there exists a section $P^*|Q$ such that 
\[ 
M'\cdot \bm\gamma^- = -[P], \ M' \cdot  \bm\gamma^+ = [Q] \ \text{ or } \
M'\cdot \bm\gamma^- = [Q^\vee], \ M' \cdot  \bm\gamma^+ = - [P^\vee].   \qedhere
\]
\end{proof}

Let $\xi, \eta \in \hat{\mathbb R}$ be two distinct points on the boundary of $\mathbb H$ and $U$ be the associated doubly-infinite $\mathbf H_4$-sequence of the oriented geodesic $\bm\gamma$ with $\bm\gamma^- = \xi$, $\bm\gamma^+=\eta$. 
For any section $P^*|Q$ of $U$, 
there exists $M \in \mathbf H_4$ such that 
$$\tilde{\bm\gamma}^- = M \cdot \xi = - [P], \quad \tilde{\bm\gamma}^+ = M \cdot \eta = [Q] \quad \text{ for } \ 
\tilde{\bm\gamma} = M(\bm\gamma).$$
Since 
$$ SM\cdot \xi = [P^\vee], \quad SM\cdot \eta = - [Q^\vee] \quad  \text{ and } \quad SM \in \mathbf H_4,$$ 
we have 
$$
\sup_{M \in \mathbf H_4} \left| M\cdot \xi - M \cdot \eta \right| \ge \sup_{P^*|Q} \left( \max \left\{ [Q] + [P], [Q^\vee] + [P^\vee] \right\} \right)\ge 2,
$$
where $P^*|Q$ runs over all sections of $U$ and the second inequality is from 
$$\left( [P] + [Q] \right)\left( [P^\vee] + [Q^\vee]\right) = 2 + \frac{[P]}{[Q]} + \frac{[Q]}{[P]} \ge 4.$$
Therefore, Lemma~\ref{lem:min2} implies that
$$
\sup_{M \in \mathbf H_4} \left| M\cdot \xi - M \cdot \eta \right| =  \sup_{P^*|Q} \left( \max \left\{ [Q] + [P], [Q^\vee] + [P^\vee] \right\} \right).$$
Let
\begin{equation*}
L(P^{*}|Q): = [P] + [Q].
\end{equation*}
Then we have Perron's formula for the Hecke group $\mathbf H_4$ as follows.

\begin{theorem}\label{P1}
Let $U$ be a doubly-infinite $\mathbf H_4$-sequence. 
We define $\mathcal M(U)$ by the maximum of two supremum values as follows:
$$ \mathcal M(U) : =  
\sup_{P^{*}|Q} \max \left\{ L(P^{*}|Q), L((P^\vee)^{*}|Q^\vee) \right\},$$
where $P^*|Q$ runs over all sections of $U$.
The Markoff spectrum is the set of 
$\mathcal M(U)$ as $U$ runs through all of doubly-infinite $\mathbf H_4$-sequences 
$$\mathscr M(\mathbf H_4) = \{ \mathcal M(U) \in \mathbb R \, | \,  U \text{ is a doubly-infinite $\mathbf H_4$-sequence} \}.
$$
\end{theorem}

\begin{theorem}\label{P2}
Let $U$ be a doubly-infinite $\mathbf H_4$-sequence. 
We define $\mathcal L(U)$ by the maximum of two limit superior values as follows:
$$ \mathcal L(U) : =  \limsup_{P^{*}|Q} \max \left\{ L(P^{*}|Q), L((P^\vee)^{*}|Q^\vee) \right\},$$
where $P^*|Q$ runs over all sections of $U$.
The Lagrange spectrum is the set of 
$\mathcal L(U)$ as $U$ runs through all of doubly-infinite $\mathbf H_4$-sequences
\begin{align*}
\mathscr L(\mathbf H_4) &= \{ \mathcal L(U) \in \mathbb R \, | \,  U \text{ is a doubly-infinite $\mathbf H_4$-sequence} \} . 
\end{align*}
\end{theorem}

For a finite sequence $W$, we denote $k$ repeated sequence $W\cdots W$ by $W^{k}$.
We also denote an infinite sequence with period $W$ and a doubly infinite sequence with period $W$ by $W^{\infty}$ and ${}^{\infty}W^{\infty}$.

\begin{example}\label{examp}
The associated doubly infinite $\mathbf H_4$-sequences of the three closed geodesics in Figure~\ref{Fig:mod2}
are $U_1 = {}^\infty2^\infty$ (left), $ U_2 = {}^\infty(13)^\infty$ (center) and $U_3= {}^\infty(123)^\infty$ (right).
From \eqref{period_e1} and \eqref{period_e2}, we check 
\begin{align*}
\mathcal M(U_1) &= L( \dots 2 2 | 22 \dots ) = 2 [2^\infty] = 2, \\
\mathcal M(U_2) &= L( \dots 3 1 3 1 | 3 1 3 1 \dots ) 
= \sqrt 2 + 2 [(13)^\infty] =  \sqrt 6, \\
\mathcal M(U_3) &= L( \dots 1 2 3 | 1 2 3 1 2 3 \dots ) 
= [(123)^\infty] + \frac{1}{[(123)^\infty]}
= \frac{2\sqrt{17}}{3}.
\end{align*}
\end{example}

Hereafter, commas in the $\mathbf H_4$-sequences may occasionally be omitted for simplicity of notation.


\section{Closedness of the Markoff spectrum}\label{Sec:closedness}

We prove Theorem~\ref{thm-1}. 
First we note that given the discrete topology on $\{1,2,3\}$, 
the product space $\{ 1, 2, 3 \}^{\mathbb Z}$ 
is compact due to Tychonoff's theorem.

\begin{lemma}\label{MarkoffTrans}
Let $U$ be a doubly-infinite $\mathbf H_4$-sequence.
If $\mathcal M(U)$ is finite, then there exists a doubly-infinite $\mathbf H_4$-sequence $\tilde U$ with a section $P^* | Q$ such that $\mathcal M(U)= \mathcal M(\tilde U) = L( P^* | Q )$.
\end{lemma}

\begin{proof}
By Theorem~\ref{P1}, there exists a sequence of sections  $\{P_n^*|Q_n\}_{n\in\mathbb N}$ of $U$ or $U^\vee$, say $U$,  
such that $\lim\limits_{n\to\infty} L(P_n^*|Q_n) = \mathcal M(U)$.
Since the product space 
$\{ 1,2,3\}^{\mathbb Z}$ is compact, there exists a subsequence $\{P_{n_k}^*|Q_{n_k}\}_{k \in \mathbb N}$ which converges to a section $P^* | Q$ of a doubly-infinite $\mathbf H_4$-sequence $\tilde U$.
By the continuity of $L$, we have $L(P^* | Q) = \mathcal M(U) \le \mathcal M(\tilde U)$.

If $\tilde P^*|\tilde Q$ is another section of $\tilde U$, then $\tilde P^*|\tilde Q$ is a limit of $\{\tilde P_{n_k}^*|\tilde Q_{n_k}\}_{k \in \mathbb N}$, which is a shifted subsequence of $\{P_{n_k}^*|Q_{n_k}\}$.
Thus $L(\tilde P^*|\tilde Q) \le \mathcal M(U)$, which implies that 
$\mathcal M(\tilde U) \le \mathcal M(U)$.
\end{proof}

\begin{proof}[Proof of Theorem~\ref{thm-1}]
We first show that the Markoff spectrum $\mathscr M(\mathbf{H}_4)$ is closed.
Choose a convergent sequence $\{m_n\}_{n\in \mathbb N}$ in $\mathscr{M}(\mathbf{H}_4)$. 
By Lemma \ref{MarkoffTrans}, there exists a sequence of doubly-infinite $\mathbf H_4$-sequences $\{U_n\}_{n \in \mathbb N}$ with a sequence of sections of $\{P_n^* | Q_n\}_{n \in \mathbb N}$ 
such that $m_n = \mathcal M(U_n)= L( P_n^* | Q_n )$ for all $n \in\mathbb{N}$. 
By the compactness of $\{ 1,2,3\}^{\mathbb Z}$, we have a converging subsequence $\{P_{n_k}^* | Q_{n_k}\}_{k \in \mathbb N}$ to the limit $P^*|Q$ which is a section of a doubly-infinite $\mathbf H_4$-sequence $U$.
By the continuity of $L$, we have $\lim\limits_{n\to\infty} m_n = L(P^*|Q)$, thus
$\lim\limits_{n\to\infty} m_n \le \mathcal M(U)$.

Let $\tilde P^*|\tilde Q$ be another section of $U$. 
Then $\tilde P^*|\tilde Q$ is a limit of finite shifts of the subsequence $\{P_{n_k}^*|Q_{n_k}\}_{k \in \mathbb N}$.
Therefore, $L(\tilde P^*|\tilde Q) \le \lim\limits_{n \to \infty} \mathcal M(U_n)$ and 
$\mathcal M(U) \le \lim\limits_{n \to \infty} m_n$.
Hence, $\mathcal M(U) = \lim\limits_{n \to \infty} m_n$ and we conclude that the Markoff spectrum is closed.

Now we show that $\mathscr{L}(\mathbf{H}_4) \subset \mathscr{M}(\mathbf{H}_4)$.
By Theorem~\ref{P2}, for a doubly-infinite $\mathbf H_4$-sequence $U$, there exists a sequence of sections $\{P_n^* | Q_n\}_{n \in \mathbb N}$ of $U$ or $U^\vee$, say $U$, 
such that $\mathcal{L}(U) = \lim\limits_{n \to \infty} L (P_n^* | Q_n)$.
Since the product space $\{ 1,2,3\}^{\mathbb Z}$ is compact, there exists a subsequence $\{P_{n_k}^* | Q_{n_k}\}_{k \in \mathbb N}$ which converges to an element $P^* | Q \in \{ 1, 2, 3 \}^{\mathbb Z}$,
which is a section of a doubly-infinite sequence $\tilde U$.
By the continuity of $L$, we deduce that 
$\mathcal L(U) \le \mathcal{M}(\tilde U)$.

For another section $\tilde P^*|\tilde Q$ of $\tilde U$,   
we have $L(\tilde P^*|\tilde Q) \le \mathcal L(U)$ since 
$\tilde P^*|\tilde Q$ is a limit of a sequence of sections of $U$.
Therefore, $\mathcal M(\tilde U) \le \mathcal L(U)$. 
Hence, $\mathcal L(U) = \mathcal M(\tilde U) \in \mathscr M(\mathbf{H}_4)$.
\end{proof}

\section{Hausdorff dimension of the Lagrange spectrum}\label{Sec:dim}

In this section, we prove Theorem~\ref{thm0}.
By \eqref{eq:M1},
for each $\mathbf H_4$-sequence $P$, $Q$, we have 
$$
\left| [1 P] - [1 Q] \right| \le \left| [P] - [Q] \right|, \quad
\left| [2 P] - [2 Q] \right| \le \frac{\left| [P] - [Q]\right|}{\sqrt 2}.
$$
Assume that $\varepsilon >0$ is given in this section.
Then we choose $m \ge 1$ such that 
$$
[ (1 2^{m} 3)^\infty ] - [ 1 2^\infty] 
\le \frac{[3(12^m3)^\infty] - [2^\infty]}{(\sqrt 2)^m}
< \varepsilon.
$$
We have for any $\mathbf H_4$-sequence $P$
\begin{equation}\label{ebound2}
\begin{split}
[ 3 2^{m+1} 1 P ] + 
[ (1 2^{m} 3)^\infty ] 
&< [3 2^\infty] + 
[1 2^\infty] + \varepsilon = 
2\sqrt 2 + \varepsilon.
\end{split}
\end{equation}
Let $A = 3 2^{m+1} 1$, $B = 32^{m}1$. 
Define 
\begin{equation*}
\Sigma := \{ P \in \{ 1,2,3 \}^{\mathbb N} 
\, | \, P = B^{m_1} A^{n_1} B^{m_2} A^{n_2} \cdots, \  n_i, m_i \in \{1,2\} \, \text{ for all $i$} \, \}.
\end{equation*}

\begin{lemma}\label{dimension}
Let $\mathcal F = \{ [P] \in \mathbb R \, | \, P \in \Sigma \}$. 
Then we have 
$$\dim_H ( \mathcal F ) > 0.$$
\end{lemma}

\begin{proof}
Let
$$
\alpha := [(B^2 A)^\infty ], \quad  \beta := [(B A^2)^\infty]. $$
Then for each $P \in \Sigma$, we have
$$
\alpha \le [P] \le \beta.$$
Let 
$$
M_A = N_3 N_2^{m+1} N_1, \qquad M_B = N_3 N_2^{m} N_1.
$$
Define $f_i : [ \alpha, \beta ] \to [ \alpha, \beta ]$ to be
$$f_1(x) = M_B^2M_A \cdot x, \ 
f_2(x) = M_B^2M_A^2 \cdot x, \
f_3(x) = M_BM_A \cdot x, \
f_4(x) = M_BM_A^2 \cdot x.$$
Then $\{f_1, f_2, f_3, f_4\}$ is a family of contracting functions, which is called  an iterated function system (see e.g. \cite{Fal03})
satisfying 
$$\mathcal F = f_1(\mathcal F) \cup f_2(\mathcal F) \cup f_3(\mathcal F) \cup f_4(\mathcal F),
\quad f_i(\mathcal F) \cap f_j(\mathcal F) = \emptyset \ \text{ for } i\ne j.$$
Using \eqref{eq:M1}, we check that there are $c_i >0$ for each $i = 1,2,3,4$ such that $|f_i(x) - f_i(y) | \ge c_i |x-y|$ for $x,y \in [ \alpha, \beta ]$ since all element of the matrices $M_B^2M_A$, $M_B^2M_A^2$, $M_BM_A$, $M_BM_A^2$ are positive.
By \cite{Fal03}*{Proposition 9.7},
we conclude that 
$$
\dim_H( \mathcal F ) \ge s, 
$$
where $s >0$ is the constant satisfying
\begin{equation*}
c_1^s + c_2^s + c_3^s + c_4^s = 1. \qedhere
\end{equation*}
\end{proof}

Choose 
$$R = B^{m_1} A^{n_1} B^{m_2} A^{n_2} \cdots \in \Sigma, \quad n_i, m_i \in \{1,2\}$$
and let
\begin{align*}
W^R_k &:= B^{m_1} A^{n_1} B^{m_2} A^{n_2} B^{m_3} \cdots B^{m_k} A^{n_k}, \\
U_R &:= {}^\infty B W^R_1 B^{2} A^3 W^R_{2} B^{3} A^3 W^R_{3} B^{4} 
\cdots B^{k} A^3 W^R_{k} B^{k+1} \cdots.
\end{align*}

\begin{lemma}\label{sections}
We have
$$
\mathcal L( U_R) = [(B^\vee)^\infty] + [A^3 R] = \frac{1}{[B^\infty]} + [A^3 R].
$$
\end{lemma}

\begin{proof}
Let $P^*32^k | 2^\ell 1 Q$ be a section of $U_R$ for some $k, \ell \ge 0$. 
Then we have for $k \ge 1, \ell \ge 0$  
\begin{equation*}
L(P^*32^k | 2^\ell 1 Q) = [2^k 3 P] + [2^\ell 1 Q] < [23^\infty] + [2^\infty] = \sqrt 2 + 1 
\end{equation*}
and for $k = 0$
\begin{equation*}
L(P^*3 2^k| 2^\ell 1 Q) =  L(P^*| 32^\ell 1 Q).
\end{equation*}
On the other way, if $P^*32^k 1 | 3 2^\ell 1 Q$ is a section of $U$, then
\begin{equation*}
L(P^*32^k 1 | 3 2^\ell 1 Q) = [1 2^k 3 P] + [3 2^\ell 1 Q] > [1 2^\infty] + [3 2 1^\infty] = \frac{5}{\sqrt 2} -1 > \sqrt 2 + 1. 
\end{equation*}
Therefore, 
we have
\begin{equation*}
\mathcal L( U_R) = \limsup_{P^*|Q} \max \left\{ L( P^*|Q), L( (P^\vee)^*|Q^\vee) \right\} 
\end{equation*}
where $P^*|Q$ runs over all sections of $U_R$ such that 
$P^*|Q = \tilde P^* A|A \tilde Q$, $\tilde P^* A|B \tilde Q$,  $\tilde P^* B|A \tilde Q$, or $\tilde P^* B|B \tilde Q$ for some $\tilde P$ and $\tilde Q$.
Using the facts that 
$[A P] > [B Q]$ for any infinite sequences $P$, $Q$
and that $W_k^R$ does not contain $A^3$,
we conclude that 
\begin{align*}
\mathcal L( U_R) 
&= \limsup_{k \to \infty} L( \cdots B_{k-1} A^3 W^R_{k-1} B^{k} | A^3 W^R_k B^{k+1} A^3 W^R_{k+1} \cdots ) \\
&= \limsup_{k \to \infty} \left( \frac{1}{[B^{k}W^R_{k-1} A^3B_{k-1} \cdots]} + [A^3 W^R_k B^{k+1} A^3 W^R_{k+1} \cdots ] \right) \\
&= L( {}^\infty B | A^3 R) = [(B^\vee)^\infty] + [A^3 R] = \frac{1}{[B^\infty]} + [A^3 R]. \qedhere
\end{align*}
\end{proof}

Let  
$$
\mathcal H : = 
\{ \mathcal L( U_R) \, | \, R \in \Sigma \}.
$$
Then, Lemma~\ref{sections} and \eqref{ebound2} yield that 
\begin{equation}
\mathcal{H} = \left\{ \frac{1}{[B^\infty]} + [A^3 R] \, | \,  R \in \Sigma \right\}
\subset \mathscr L(\mathbf H_4) \cap  (0, 2\sqrt 2 + \varepsilon).
\end{equation}
Since all elements of the matrix $M_A^3$ are positive, the map
$[R] \mapsto [A^3 R ] = M_A^3 \cdot [R]$ is a bi-Lipschitz function on the closed interval $[\alpha, \beta]$. 
Therefore, Lemma~\ref{dimension} implies that 
$\dim_H(\mathcal H) >0$ and we complete the proof of Theorem~\ref{thm0}.

\section{Gaps of the Markoff spectrum}\label{Sec:Gaps}
We investigate the gaps of $\mathscr{M}(\mathbf{H}_4)$ above the first limit point $2\sqrt 2$ in this section.
We prove Theorem~\ref{thm1} through Theorems~\ref{GapContaining232} and \ref{GapBelow232}.

We check that
\begin{equation}\label{period_e3} 
[(21)^{\infty}] = \frac{\sqrt 2}{\sqrt{7} - 1}, \quad
[(2131)^\infty] = \frac{\sqrt{119} + 3}{11\sqrt 2} 
\end{equation}
and let 
\begin{align*}
m_0 &:= \mathcal{M} \left( {}^\infty(3132) 123 (2131)^\infty \right) \\
&= L ({}^\infty(3132) 12 \,| \, 3 (2131)^\infty) = \REP = 3.181 \dots. 
\end{align*}

\begin{lemma}\label{lem:gap}
Let $\mathcal M(U) \le m_0$.
Then $U$ satisfies one of the followings: 
\begin{enumerate}[label=\upshape(\roman*), leftmargin=*, widest=iii]
\item $U = {}^\infty(1232)^\infty$, or 
\item $U$ or $U^\vee = {}^\infty(3132) 123 (2131)^\infty$, or 
\item $U$ does not contain 11, 33, 212, 232.
\end{enumerate}
\end{lemma}

\begin{proof}
First, if $U$ or $U^\vee$, say $U$, contains 333, then 
$$\mathcal{M}(U) \ge L(P^{*}|333Q) = [P] + [333Q] = [P] + [Q] + 3\sqrt 2 \ge 3 \sqrt 2 >m_0$$ 
for some infinite $\mathbf H_4$-sequences $P,Q$ with $U=P^{*}333Q$.
Therefore, $U$ and $U^\vee$ do not contain 333 nor 111.

Next, assume that $U$ or $U^\vee$, say $U$, contains 33.
Let $U=P^{*}33Q$ for some infinite $\mathbf H_4$-sequences $P,Q$ starting with 1 or 2. 
Then, by \eqref{period_e2}, we have
\begin{align*}
\mathcal{M}(U) \ge L(P^{*}|33Q) &= [Q] + [P] + 2\sqrt 2
\ge [(112)^\infty] + [(112)^\infty] + 2\sqrt 2 \\
&= \frac{2}{\sqrt{7} +\sqrt2} +2\sqrt 2 > m_0.
\end{align*}
Hence, $U$ and $U^\vee$ do not contain 33 nor 11.

We claim that $U$ and $U^*$ do not contain 2322 or 2323. 
Let $U=P^{*}232Q$ for some infinite $\mathbf H_4$-sequences $P$,$Q$ with $Q$ beginning with 2 or 3.
Then,  by 
\eqref{period_e3}, 
\begin{align*}
\mathcal{M}(U) \ge L(P^{*} 2|32 Q)
&= [2P] + [2Q] + \sqrt 2 
\ge [2(12)^\infty] + [2 (2 1)^\infty] + \sqrt 2 \\
&= \frac{\sqrt 2}{\sqrt{7} - 1} + \frac{\sqrt{7} + 1}{\sqrt{14}} + \sqrt 2 > m_0.
\end{align*}
Therefore, $U$ does not contain 11, 33, 2323, 2121, 3232, 1212, 2232, 2212, 2322, 2122. Thus any infinite $\mathbf H_4$-sequence $P$ of $U$ or $U^*$ satisfies that 
\begin{equation}\label{lembound}
\left[ (1213)^\infty \right] \le [P] \le [ (3231)^\infty ].
\end{equation}

Suppose that $U \ne {}^\infty(1232)^\infty$ and $U$ contains 232. 
%
%
%
Then $U$ contains (a) 3123213, or (b) 2123212 or, (c) 3123212 or 2123213, say 3123212.

(a) If $U$ contains 3123213, then $U=P^{*}3123213Q$ for some infinite $\mathbf H_4$-sequences $P$,$Q$ and by 
\eqref{lembound} and \eqref{period_e3} we have
\begin{align*}
\mathcal{M}(U) &\ge L(P^{*} 312|3213 Q)
= [213P] + [213Q] + \sqrt 2  \\
&\ge [(2131)^\infty] + [(2131)^\infty] + \sqrt 2
= 2 \frac{\sqrt{119} + 3}{11\sqrt 2} + \sqrt 2 > m_0.
\end{align*}
Hence, $U$ does not contain 3123213 nor 1321231.

(b) Assume that $U$ contains 2123212. 
Since $U \neq {}^\infty(1232)^\infty$, 
there exists an infinite $\mathbf H_4$-sequence $P$ not beginning with 32, for which  
$U$, $U^\vee$, $U^*$, or $(U^\vee)^*$, say $U^\vee$,
is $P^{*}2123212Q$ 
for some $Q$.
Hence, by \eqref{lembound}, we have 
\begin{align*}
\mathcal{M}(U) 
&\ge L((P^\vee)^{*}2|321232Q^\vee)
= [2P^\vee] + [21232Q^\vee] + \sqrt 2 \\
& > [2 (1312)^\infty] +[212(3132)^\infty] + \sqrt 2\\
& =\mathcal{M}({}^{\infty}(3132)123 (2131)^\infty)=m_0.
\end{align*}

(c) If $U$ contains 2123213, but does not contain 2123212 nor 2321232, then 
we have $U=P^{*}21232Q$ for some infinite $\mathbf H_4$-sequences $P,Q$, where $P$ does not begin with 32 and $Q$ does not start with 12.
Thus, \eqref{lembound} implies  
$$[P^\vee], \, [Q] \ge [13(1213)^\infty] = [(1312)^\infty].
$$
By the elementary calculus, we check that $[2P]+[212P^\vee]$ is an increasing function of $[P]$ on the interval
$\left( [(1312)^\infty], [(3132)^\infty] \right)$.
Therefore, 
we have 
\begin{align*}
\mathcal{M}(U) 
&\ge \frac{1}{2}(L(P^{*}212|32Q)+ L((P^\vee)^{*}23|212Q^\vee)) \\
&= \frac 12 \left( [212P] + [2P^\vee] \right) + \frac 12 \left( [2Q] + [212Q^\vee] \right) + \sqrt 2 \\
&\ge \frac{[2(1312)^\infty]+[212(3132)^\infty]}{2} +\frac{[2(1312)^\infty]+[212(3132)^\infty]}{2} + \sqrt 2 \\
&=[(2131)^\infty]+[21(2313)^\infty] + \sqrt 2 
=\mathcal{M}({}^{\infty}(3132)123(2131)^\infty) =m_0. 
\end{align*} 
Moreover, if the equality holds, then $U = {}^{\infty}(3132)123(2131)^\infty$ or $U^\vee = {}^{\infty}(3132)123(2131)^\infty$.
\end{proof}

\begin{theorem}\label{GapContaining232}
The interval 
$$
\left( \sqrt{10}, \REP \right)=(3.162\dots, 3.181\ldots)
$$ 
is a maximal gap in $\mathscr{M}(\mathbf{H}_4)$.
Two boundary points of the interval correspond to $\mathcal{M} 
\left({}^\infty(1232)^\infty \right)=\sqrt{10}$ and $\mathcal{M}\left( U \right)=\REP$ for $U = {}^\infty(3132)123(2131)^\infty$.
Moreover, $\REP$ is a limit point of $\mathscr{M}(\mathbf{H}_4)$.
\end{theorem}

\begin{proof}
By direct calculation, we check that 
	$$\mathcal{M}({}^\infty (1232)^\infty)=\sqrt{10} \ \text{ and } \ \mathcal{M}( {}^\infty(3132)123(2131)^\infty )=\REP =: m_0.$$

Suppose that $U$ is a doubly infinite $\mathbf H_4$-sequence such that $\mathcal{M}(U) < m_{0}$ and $U \neq {}^\infty (1232)^\infty$.
Then, 
by Lemma~\ref{lem:gap}, $U$ does not contain 11, 33, 212, 232.
If $a, b \in \{ 1,2 \}$, then  
\begin{equation*}
L(P^{*}a \, |\, bQ) = [aP] + [bQ] \le [2P]+[2Q] \le \sqrt 2 + \sqrt 2 <\sqrt{10}
\end{equation*}
for any infinite $\mathbf H_4$-sequences $P,Q$.
Moreover, if $P^{*} 1 | 3a Q$ is a section of $U$ for $a \in \{ 1,2\}$, then 
\begin{align*}
L(P^{*} 1 \, | \, 3a Q)
&= [1P] + [3 a Q] \le
[(1323)^\infty] + [(3231)^\infty] \\ 
&=\frac{\sqrt{119}-7}{5\sqrt 2} + \frac{\sqrt{119}+7}{5\sqrt 2} = \frac{\sqrt{238}}{5}
<\sqrt{10}.
\end{align*}
Therefore, we deduce that $\mathcal M(U) < \sqrt{10}$.
Hence, $(\sqrt{10}, m_{0})$ is a maximal gap in $\mathscr{M}(\mathbf{H}_4)$.

Finally, let us show that $m_{0}$ is a limit point of $\mathscr{M}(\mathbf{H}_4)$. 
For $k\ge1$, let 
$$U_{k}:= {}^\infty(1312) 321 A_k 123 (2131)^\infty,
\ \text{ where } \ A_{k}:=(2313)^k 2 = 2(3132)^k.
$$ 
Let $P_{k}^* | Q_k$ be a section of the doubly infinite sequence $U_k$.
Then there exists a section $P^* | Q$ of the doubly infinite sequence $U$ such that at least the first $4k+2$ digits of $P$, $P_k$ and those of $Q$ and $Q_k$ are identical. 
Therefore, we have 
$\lim\limits_{k\to\infty}{\mathcal{M}(U_{k})}=\mathcal{M}(U) = m_{0}$. 
By Lemma~\ref{lem:gap}, we have $\mathcal{M}(U_{k}) >  m_{0}$ for all $k$. 
Hence, $m_{0}$ is a limit point of $\mathscr{M}(\mathbf{H}_4)$.
\end{proof}

\begin{theorem}\label{GapBelow232}
The interval 
\[ \left(\frac{\sqrt{238}}{5}, \sqrt{10} \right) = (3.085\dots, 3.162\dots )\]
is a maximal gap in $\mathscr{M}(\mathbf{H}_4)$. 
The lower end point satisfies $\mathcal{M}({}^{\infty}(1312)^{\infty})=\frac{\sqrt{238}}{5}$. 
\end{theorem}

\begin{proof}
Suppose that $\mathcal{M}(U) < \sqrt{10} = \mathcal{M} \left({}^\infty(1232)^\infty \right)$.
By Lemma~\ref{lem:gap}, $U$ does not contain 33,11,212,232.
Therefore, for any infinite $\mathbf H_4$-sequence $R$ appearing in $U$, if $R$ does not start with 32 or 12, then 
$$[(1312)^\infty] \le [R] \le [(3132)^\infty].$$ 
If $U=P^{*}2Q$ for infinite $\mathbf H_4$-sequences $P,Q$, then both $P$, $Q$ do not start with 32. 
Therefore, 
\begin{equation*}
L(P^{*}|2Q) = [P] + [2Q] \le [(3132)^\infty]+[2 (3132)^\infty] = \mathcal{M}({}^\infty(3132)^\infty).
\end{equation*}
If $U=P^{*}13Q$, then
$[Q]\le[(2313)^\infty]$ and $[31P]\le[(3132)^\infty]$. 
Thus,
\begin{align*}
L(P^{*} 1 |3 Q) 
&= L(P^{*}13 | Q) 
= [31P] + [Q] \\
&\le [(2313)^\infty] + [(3132)^\infty] =\mathcal{M}({}^\infty(3132)^\infty).
\end{align*}
Therefore, for any section $P^*|Q$ of $U$, we have 
$L(P^*|Q) \le \mathcal{M}({}^\infty(3132)^\infty) =\frac{\sqrt{238}}{5}$.
Thus, by Theorem~\ref{P1}, the interval $\big(\frac{\sqrt{238}}{5}, \sqrt{10}\big)$ is a maximal gap in $\mathscr{M}(\mathbf{H}_4)$.
\end{proof}

\section{A bound of Hall's Ray}\label{Sec:Halslray}

In this section, we give the bound of Hall's ray (Theorem~\ref{thm2}).
Let
$$\mathcal{K} =\{\,[d_1,d_2, \dots] \;| \; d_1 \ne 3, \ d_kd_{k+1}d_{k+2} \ne 111 \text{ nor } 333 \ \text{ for all } k \ge 1\}
$$
and let $R= (332)^\infty$.
We note that $[R] 
= \sqrt 7 + \sqrt 2$ and 
$$
\min \mathcal{K} = [R^\vee]= \frac{\sqrt{7}-\sqrt 2}{5}, \quad \max \mathcal{K} = [2R]=\sqrt{7}-\sqrt 2.
$$
Let $\mathcal{E}(c_1, \dots, c_n)$ be the smallest closed interval containing $\{ [d_1, d_2, \dots] \in \mathcal{K} \, | \, d_1 = c_1, 
\dots, d_n = c_n \}$
and $\mathcal E = \big[ [R^\vee], [2R] \big] =\big[\frac{\sqrt{7}-\sqrt 2}{5}, \sqrt{7}-\sqrt 2\big].$
Then we have 
\begin{equation*}
\mathcal{E}(c_1, \dots, c_n) 
= \begin{cases}
\big[ \,[c_1 \cdots c_{n} 2 R^\vee], \, [c_1 \cdots c_{n} R] \, \big],  &\text{ if } \ 
c_{n-1}c_n = 11, \\
\big[\, [c_1 \cdots c_{n} 12 R^\vee],\, [c_1 \cdots c_{n} R]\, \big], &\text{ if } \ 
c_{n-1} \ne 1, c_n = 1, \\
\big[\, [c_1 \cdots c_n R^\vee], \, [c_1 \cdots c_n R]\,\big], &\text{ if } \ 
c_n = 2, \\
\big[\,[c_1 \cdots c_{n} R^\vee], \, [c_1 \cdots c_{n} 32 R]\, \big], &\text{ if } \ 
c_{n-1} \ne 3, c_n = 3, \\
\big[ \,[c_1 \cdots c_{n} R^\vee], \, [c_1 \cdots c_{n} 2 R] \, \big], &\text{ if } \ 
c_{n-1}c_n = 33.
\end{cases}
\end{equation*}
We also define 
$\mathcal{E}_*(c_1, \dots, c_n)$ be the smallest closed interval containing 
$$\{ [d_1, d_2, \dots] \in \mathcal{K} \, | \, d_1 = c_1, d_2 = c_2, \dots, d_n = c_n, d_{n+1} \neq 3 \}.$$

First, let us verify that $\mathcal{K}$ can be obtained by applying the Cantor dissection process.
In the dissection process, each type of interval is divided by the following rules:

(i) The interval $\mathcal{E}(c_1, \dots, c_n)$ is divided into the union of two intervals 
$$
\begin{cases}
\mathcal{E}(c_1, \dots, c_n, 2) \cup \mathcal{E}(c_1, \dots, c_n, 3), &\text{ if } c_{n-1}c_n = 11, \\
\mathcal{E}(c_1, \dots, c_n, 1) \cup \mathcal{E}(c_1, \dots, c_n, 2), &\text{ if } c_{n-1}c_n = 33, \\
\mathcal{E}_*(c_1, \dots, c_n) \cup \mathcal{E}(c_1, \dots, c_n, 3), &\text{ otherwise}. \\
\end{cases}
$$

(ii) The interval $\mathcal{E}_*(c_1, \dots, c_n)$ with $c_{n-1} c_n \ne 11$ nor $33$, is divided into the union of two intervals 
$$\mathcal{E}(c_1, \dots, c_n, 1) \cup \mathcal{E}(c_1, \dots, c_n, 2).$$

Each type of interval is subdivided into two intervals. 
Starting from $\mathcal{E}$, we continued the dissection process according to the above rules. 
Thus, we obtain the Cantor set $\mathcal K$.

\begin{lemma}\label{IntervalCond}
Let $\mathcal{I}$ be a closed interval of $\mathcal{E}(c_1, \dots, c_n)$ or $\mathcal{E}_*(c_1, \dots, c_n)$.
In the Cantor dissection process,
we have closed intervals $\mathcal{I}_{1}, \mathcal{I}_{2}$ in $\mathcal{I}$ satisfying $\mathcal{I} \setminus \mathcal{J} = \mathcal{I}_1 \cup \mathcal{I}_2$ for an open interval $\mathcal{J}$. 
Then
$$|\mathcal{I}_{i}| \ge |\mathcal{J}| \quad \textrm{for}\quad i=1,2. $$
\end{lemma}

\begin{proof}
Let
$$M = N_{c_1} N_{c_2} \cdots N_{c_n} = \begin{pmatrix} p & r \\ q & s \end{pmatrix}. $$
Then we note that 
from \eqref{eq:M1}, we have
$$
[c_1\dots c_n  P] - [c_1 \dots c_n Q] 
= M \cdot [P] - M \cdot [Q] = \frac{[P] - [Q]}{(q[P]+s)(q[Q]+s)}.
$$
Using \eqref{invN}, we have
$$
\frac{s}{q} = - M^{-1} \cdot \infty = - N^{-1}_{c_n} \cdots N^{-1}_{c_1} \cdot \infty = N_{c_n} \cdots N_{c_1} \cdot \infty
= [c_n \cdots c_1 3^\infty], $$
we have 
\begin{align}
\frac sq
&\le [11 3^\infty] = \frac{1}{2\sqrt 2}, \quad \text{ if } \ c_{n-1}c_n = 11, \label{matrix_b1} \\
\frac sq &> [11 3^\infty] = \frac{1}{2\sqrt 2}, \quad \text{ if } \  c_{n-1}c_n \ne 11. \label{matrix_b2}
\end{align}


Let $\mathcal{I} = \mathcal{E}(c_1, \dots, c_n)$ with $c_{n-1} c_n \ne 11$ nor $33$.
Then we have
$$\mathcal{I}_1 = \mathcal{E}_*(c_1, \dots, c_n), \quad \mathcal{I}_2 = \mathcal{E}(c_1, \dots, c_n, 3).$$
Therefore, we have
\begin{align*}
|\mathcal{J}| &= [c_1 \dots c_{n} 3 R^\vee] - [c_1 \dots c_{n} 2 R]  
=\frac{[3R^\vee]-[2R]}{( q[3R^\vee]+s)(q[2R]+s)}, \\
|\mathcal{I}_1| &\ge [c_1 \dots c_{n} 2 R] - [c_1 \dots c_{n} 12R^\vee] 
= \frac{[2R]-[12R^\vee]}{(q[2R]+s)(q[12R^\vee]+s)}, \\
|\mathcal{I}_2| &\ge [c_1 \dots c_{n}32R] - [c_1 \dots c_{n} 3 R^\vee] 
= \frac{[32R]-[3R^\vee]}{(q[32R]+s)(q[3R^\vee]+s)},
\end{align*}
thus
\begin{align*}
\frac{|\mathcal{J}|}{|\mathcal{I}_1|} &\le \frac{(q[12R^\vee]+s)([3R^\vee]-[2R])}{( q[3R^\vee]+s)([2R]-[12R^\vee])} < \frac{[3R^\vee]-[2R]}{[2R]-[12R^\vee]} 
= 0.5025 \dots < 1, \\
\frac{|\mathcal{J}|}{|\mathcal{I}_2|} &\le \frac{(q[32R]+s)([3R^\vee]-[2R])}{( q[2R]+s)([32R]-[3R^\vee])} < \frac{[32R]([3R^\vee]-[2R])}{ [2R]([32R]-[3R^\vee])} \\
&= 0.9354\dots < 1.
\end{align*}

Let $\mathcal{I} = \mathcal{E}(c_1, \dots, c_n)$ with $c_{n-1}c_n=11$.
Then we have
$$\mathcal{I}_1 = \mathcal{E}(c_1, \dots, c_n, 2), \quad \mathcal{I}_2 = \mathcal{E}(c_1, \dots, c_n, 3).$$
Therefore, 
\begin{align*}
|\mathcal{J}| &= [c_1\dots c_{n} 3 R^\vee] - [c_1 \dots c_{n} 2 R] 
=\frac{[3R^\vee]-[2R]}{(q[3R^\vee]+s)(q[2R]+s)}, \\
|\mathcal{I}_1| &= [c_1 \dots c_{n} 2 R] - [c_1 \dots c_{n}2R^\vee] 
=\frac{[2R]-[2R^\vee]}{(q[2R]+s)(q[2R^\vee]+s)}, \\
|\mathcal{I}_2| &= [c_1 \dots c_{n}R] - [c_1 \dots c_{n} 3 R^\vee] 
= \frac{[R]-[3R^\vee]}{(q[R]+s)(q[3R^\vee]+s)}.
\end{align*}
By \eqref{matrix_b1}, we have 
$q \ge 2\sqrt 2 s$, thus 
\begin{align*}
\frac{|\mathcal{J}|}{|\mathcal{I}_1|} &= \dfrac{(q[2R^\vee]+s)([3R^\vee]-[2R])}{( q[3R^\vee]+s)([2R]-[2R^\vee])} \le \frac{2 \sqrt 2 [2R^\vee] + 1}{2 \sqrt 2 [3R^\vee] + 1} \cdot \dfrac{[3R^\vee]-[2R]}{[2R]-[2R^\vee]} \\
&=0.5917\dots <1, \\
\frac{|\mathcal{J}|}{|\mathcal{I}_2|} &= \dfrac{(q[R]+s)([3R^\vee]-[2R])}{( q[2R]+s)([R]-[3R^\vee])} \le \dfrac{[R]([3R^\vee]-[2R])}{ [2R]([R]-[3R^\vee])} = 0.5893\dots < 1.
\end{align*}

Let $\mathcal{I} = \mathcal{E}_*(c_1, \dots, c_n)$ with $c_{n-1} c_n \ne 11$ nor $33$ or $\mathcal{I} = \mathcal{E}(c_1, \dots, c_n)$ with $c_{n-1}c_n=33$.
Then we have
$$\mathcal{I}_1 = \mathcal{E}(c_1, \dots, c_n, 1), \quad \mathcal{I}_2 = \mathcal{E}(c_1, \dots, c_n, 2).$$
Therefore, we have
\begin{align*}
|\mathcal{J}| &= [c_1 \dots c_n 2 R^\vee] - [c_1 \dots c_n 1 R]  
=\dfrac{[2R^\vee]-[1R]}{(q[2R^\vee]+s)(q[1 R]+s)}, \\
|\mathcal{I}_1| &\ge [c_1 \dots c_{n} 1 R] - [c_1 \dots c_{n} 12R^\vee] 
=\dfrac{[1R]-[1 2 R^\vee]}{(q[1 R]+s)(q[12R^\vee]+s)}, \\
|\mathcal{I}_2| &= [c_1\dots c_{n}2R] - [c_1 \dots c_{n} 2 R^\vee] 
= \dfrac{[2R]-[2R^\vee]}{(q[2R]+s)(q[2R^\vee]+s)}.
\end{align*}
Using the condition that $c_{n-1} c_n \ne 11$, 
\eqref{matrix_b1} implies $q < 2\sqrt 2 s$. 
Thus,
\begin{align*}
\frac{|\mathcal{J}|}{|\mathcal{I}_1|} &\le \dfrac{(q[12R^\vee]+s)([2R^\vee]-[1R])}{( q[2R^\vee]+s)([1R]-[12R^\vee])} < 
\dfrac{[2R^\vee]-[1R]}{[1R]-[12R^\vee]} 
= 0.9354 \dots <1, \\
\frac{|\mathcal{J}|}{|\mathcal{I}_2|} &= \dfrac{(q[2R]+s)([2R^\vee]-[1R])}{( q[1R]+s)([2R]-[2R^\vee])} <  \frac{2\sqrt 2 [2R]+1}{2\sqrt 2 [1R]+1} \cdot \dfrac{[2R^\vee]-[1R]}{[2R]-[2R^\vee]} \\
&= 0.8292 \dots <1. \qedhere
\end{align*}
\end{proof}

For $\mathcal X, \mathcal Y \subset \mathbb R$, we write $\mathcal X + \mathcal Y := \{ x + y \, | \, x \in \mathcal X , y \in \mathcal Y\}$. 
\begin{lemma}(\cite{CF89}*{Chapter 4}, Lemma 3)\label{IntervalLem2}
Let $\mathcal B$ be the union of disjoint closed intervals $\mathcal A_{1}, \mathcal A_{2}, \dots, \mathcal A_r$. 
Given an open interval $\mathcal I$ in $\mathcal A_1$, let $\mathcal A_{r+1}, \mathcal A_{r+2}$ be the disjoint closed intervals such that $\mathcal A_1 \setminus \mathcal I =\mathcal A_{r+1}\cup \mathcal A_{r+2}$. 
Let $\mathcal B'$ be the union of $\mathcal A_{2}, \mathcal A_{3}, \dots,\mathcal A_{r+1}, \mathcal A_{r+2}$. 
If $|\mathcal A_{i}| \ge |\mathcal I|$ for $i=2,\dots,r+2$, then
$$\mathcal B+\mathcal B=\mathcal B'+\mathcal B'.$$
\end{lemma}

\begin{lemma}(\cite{CF89}*{Chapter 4}, Lemma 4)\label{IntervalLem3}
If $\mathcal C_{1}, \mathcal C_{2}, \dots$ is a sequence of the bounded closed sets such that $\mathcal C_i$ contains $\mathcal C_{i+1}$ for all $i\ge 1$, then
$$ \cap_{i=1}^{\infty} {\mathcal C_i} + \cap_{i=1}^{\infty} {\mathcal C_i} = \cap_{i=1}^{\infty} ({\mathcal C_{i}+\mathcal C_{i}}).$$
\end{lemma}



\begin{theorem}\label{EquivMainThm}
We have $\mathcal{K}+\mathcal{K} = \left [\frac{2\sqrt{7}-2\sqrt 2}{5}, 2\sqrt{7}-2\sqrt{2} \right]$.
\end{theorem}

\begin{proof}
Let 
$\mathcal{K}_0 := \mathcal E= \big[\frac{\sqrt{7}-\sqrt 2}{5}, \sqrt{7}-\sqrt 2\big]$
and 
construct a sequence of sets $\{\mathcal{K}_k\}_{k=0}^{\infty}$ satisfying the following four properties:
\begin{enumerate}
    \item Each $\mathcal{K}_{k}$ is closed and bounded.
    \item $\mathcal{K}_k \supset \mathcal{K}_{k+1}$ for all $k\ge0$.
    \item $\bigcap_{k=0} ^{\infty} {\mathcal{K}_k} = \mathcal{K}$.
    \item $\mathcal{K}_k +\mathcal{K}_k = \mathcal{K}_{k+1} + \mathcal{K}_{k+1}$ for $k\ge0$.
\end{enumerate}
We already verified that $\mathcal{K}$ is obtained from $\mathcal{K}_0$ by removing an infinite number of disjoint open intervals.
Now, let us arrange the set of an infinite number of the deprived open intervals in decreasing order of length and denote them by $\mathcal{J}_{0}, \mathcal{J}_{1}, \dots$.
For $k\ge0$, we set $\mathcal{K}_{k+1}=\mathcal{K}_{k}\setminus \mathcal{J}_{k}$. 
By the definition of $\mathcal{K}_k$, three properties (1), (2), (3) are satisfied.
Thus, it is enough to show that $\mathcal{K}_k +\mathcal{K}_k = \mathcal{K}_{k+1} + \mathcal{K}_{k+1}$.


Let $\mathcal{I}$ be the closed interval from which $\mathcal{J}_k$ is removed and $\mathcal{I}_{1}$, $\mathcal{I}_{2}$ be the disjoint closed intervals such that $\mathcal{I} \setminus \mathcal{J}_k =\mathcal{I}_1 \cup \mathcal{I}_2$.
By Lemma \ref{IntervalCond}, $|\mathcal{I}_{1}|,|\mathcal{I}_{2}|\ge |\mathcal{J}_{k}|$.
By the ordering of the index of $\mathcal{J}_{k-1}$ and Lemma~\ref{IntervalCond}, each closed interval in $\mathcal{K}_k$ has length greater than or equal to $|\mathcal{J}_{k-1}|$.
Hence, each closed interval in $\mathcal{K}_{k+1}$ has length equal to or greater than $|\mathcal{J}_{k}|$. 
By Lemma \ref{IntervalLem2}, $\mathcal{K}_k +\mathcal{K}_k = \mathcal{K}_{k+1} + \mathcal{K}_{k+1}$. 
Therefore, by Lemma \ref{IntervalLem3}, $\mathcal{K}+\mathcal{K}=(\cap_{i=1}^{\infty} {\mathcal{K}_i}) + (\cap_{i=1}^{\infty} {\mathcal{K}_i})=\cap_{i=1}^{\infty} ({\mathcal{K}_{i}+\mathcal{K}_{i}})=\mathcal{K}_0 +\mathcal{K}_0$.
\end{proof}

Since the length of $\mathcal{K}_0 +\mathcal{K}_0= \left[\frac{2\sqrt{7}-2\sqrt 2}{5}, 2\sqrt{7}-2\sqrt 2 \right]$ is greater than $\sqrt 2$, Theorem~\ref{EquivMainThm} implies the following corollary.

\begin{corollary}\label{MainCorInHR}
Any real number can be expressed as $\sqrt2n+[P]+[Q]$ for $n\in\mathbb{Z}$, $[P] , [Q] \in \mathcal{K}$.
\end{corollary}

Now, we obtain the bound of Hall's ray:

\begin{proof}[Proof of Theorem~\ref{thm2}]
Let $\alpha \ge 4\sqrt 2$. 
By Corollary \ref{MainCorInHR}, there exist two $\mathbf H_4$-sequences $P , Q \in \mathcal{K}$ and $n\in\mathbb{Z}$ such that $\alpha = \sqrt 2 n +[P]+[Q]$.
Since $[P], [Q] \le \sqrt 7 - \sqrt2 < \sqrt 2$, we have $n \ge 3$.
We set $P = (a_1, a_2, \dots)$ and $Q = (b_{1}, b_{2},\dots)$. 
Let $m_k$ and $\ell_k$ be increasing sequences satisfying $a_{\ell_k} \ne 3$ and $b_{m_k} \ne 3$.
Put $A_k = a_1 a_2 \dots a_{\ell_k}$ and 
$B_k = b_1 b_2 \dots b_{m_k}$.
Define a doubly infinite sequence 
$$
U = {}^\infty2 A^*_1 3^n B_1 A^*_2 3^n B_2 A^*_3 3^n B_3 A^*_4 3^n B_4 \cdots
$$
Note that $A_k$, $B_k$ do not contain 333 and the first and the last digit of $A_k$, $B_k$ are not 3.
Since $L(\tilde P^*23|32 \tilde Q) \le 4\sqrt 2$ for any section $\tilde P^*23|32\tilde Q$ of $U$, 
we have 
\begin{align*}
\mathcal{L}(U) 
&= \limsup_{k \to \infty}
L({}^\infty2 A^*_1 3^n B_1 \dots A^*_{k-1} 3^n B_{k-1} A^*_k 3^n \, | \, B_k A^*_{k+1} 3^n B_{k+2} \cdots)\\
&= n\sqrt 2 + \limsup_{k \to \infty}
 \left( [ A_k B^*_{k-1} 3^n 
 \dots B^*_1 3^n A_1 {2}^\infty ] + [ B_k A^*_{k+1} 3^n B_{k+2} \cdots ] \right)\\
&= n \sqrt 2 + [P] + [Q] = \alpha.
\end{align*}
Therefore,  
$\mathscr{L}(\mathbf{H}_4)$ contains every real number greater than or equal to $4\sqrt 2$.
\end{proof}

\subsection*{Acknowledgements}
The authors would like to thank the anonymous referee for careful reading and helpful comments.
They also thank Byungchul Cha for helpful discussions and comments, and Mauro Artigiani and Luca Marchese for drawing their attention to reference [1].
D.K. was supported by the National Research Foundation of Korea (NRF-2018R1A2B6001624,  RS-2023-00245719)
and D.S. was supported by the National Research Foundation of Korea (NRF-2020R1A2C1A01011543).

\end{document}